\documentclass[11pt]{article}

\usepackage{multirow, multicol}   
\usepackage{booktabs} 
\usepackage{array} 
\usepackage{paralist} 
\usepackage{verbatim} 
\usepackage{subfig} 

\usepackage[final]{graphicx}

\usepackage{amsmath, amssymb}
\usepackage{pictexwd,dcpic}
\usepackage{enumerate}

\newtheorem{theorem}{Theorem}
\newtheorem{corollary}[theorem]{Corollary}
\newtheorem{remark}[theorem]{Remark}
\newtheorem{lemma}[theorem]{Lemma}
\newtheorem{definition}[theorem]{Definition}
\newtheorem{proposition}[theorem]{Proposition}
\newtheorem{example}[theorem]{Example}

\newcommand{\proof}{ {\sc Proof.\quad}}
\newcommand{\pend}{ \hfill $\square$ \\}

\numberwithin{equation}{section}  
\numberwithin{figure}{section}    
\numberwithin{table}{section}     
\numberwithin{theorem}{section}

\newcommand{\of}[1]{\ensuremath{\left( #1 \right)}}

\newcommand{\abs}[1]{\ensuremath{\left| #1 \right|}}
\newcommand{\cb}[1]{\ensuremath{ \left\{ #1 \right\} }}
\newcommand{\sqb}[1]{\ensuremath{ \left[ #1 \right] }}

\newcommand{\bs}{\backslash}

\newcommand{\st}{\,|\;}

\newcommand{\vp}{\ensuremath{\varphi}}

\newcommand{\R}{\mathrm{I\negthinspace R}}
\newcommand{\OLR}{\overline{\mathrm{I\negthinspace R}}}
\newcommand{\N}{\mathrm{I\negthinspace N}}

\renewcommand{\P}{\ensuremath{\mathcal{P}}}

\newcommand{\G}{\ensuremath{\mathcal{G}}}

\newcommand{\C}{\ensuremath{C^-\setminus\cb{0}}}

\newcommand{\Min}{{\rm Min\,}}

\newcommand{\Eff}{{\rm Eff}}

\newcommand{\Limsup}{{\rm Lim}\sup}

\newcommand{\dom}{{\rm dom \,}}
\newcommand{\epi}{{\rm epi \,}}

\newcommand{\gr}{{\rm graph \,}}
\newcommand{\cl}{{\rm cl \,}}

\newcommand{\co}{{\rm co \,}}
\newcommand{\cone}{{\rm cone\,}}

\newcommand{\Int}{{\rm int\,}}

\newcommand{\isum}{{+^{\negmedspace\centerdot\,}}}

\newcommand{\idif}{{-^{\negmedspace\centerdot\,}}}

\newcommand{\lel}{\preccurlyeq}

\newcommand{\triup}{{\rm \vartriangle}}

\usepackage[
colorlinks=true, linkcolor=blue, citecolor=blue, urlcolor=blue, filecolor=blue
]{hyperref}

\begin{document}
\title{Set-optimization meets variational inequalities}

\author{
Giovanni P. Crespi\thanks{University of Valle d'Aosta,
Department of Economics and Political Sciences, Loc. Grand Chemin 73-75, 11020 Saint-Christophe, Aosta, Italy.
\href{mailto:g.crespi@univda.it}{g.crespi@univda.it}}
\and
Carola Schrage\thanks{University of Valle d'Aosta,
Department of Economics and Political Sciences, Loc. Grand Chemin 73-75, 11020 Saint-Christophe, Aosta, Italy.
\href{mailto:carolaschrage@gmail.com}{carolaschrage@gmail.com}}
}
\date{{\small \today}}
\maketitle
\begin{abstract}
We study necessary and sufficient conditions to attain solutions of set-optimization problems in terms of variational inequalities of Stampacchia and Minty type.
The notion of solution we deal with has been introduced in \cite{HeydeLoehne11}. To define the set-valued variational inequality, we introduce a set-valued directional derivative and we relate it to the Dini derivatives of a family of scalar problems.
The optimality conditions are given by Stampacchia and Minty type Variational inequalities, defined both by set-valued directional derivatives and by Dini derivatives of the scalarizations. The main results allow to obtain known variational characterizations for vector-valued optimization problems as special cases.
\end{abstract}

\section{Introduction}

Since the seminal papers by F. Giannessi (see \cite{Gian80, Gian98}), variational inequalities have been applied to obtain necessary and sufficient optimality conditions in vector optimization. In \cite{HeydeLoehne11} a new approach to study set-valued problems has been applied to have a {\em fresh look} to vector optimization. Indeed, it turns out that vector optimization can be treated as a special case of set-valued optimization.
The aim of this paper is to provide some variational characterization of (convex) set-valued optimization. Following the approach known as set-optimization we introduce set-valued variational inequalities, both of Stampacchia and Minty type, by means of Dini-type derivatives (see e.g. \cite{HamelSchrage13PJO}). Under suitable assumptions (e.g. lower semicontinuity type assumptions), we can prove equivalence between solutions of the variational inequalities and solutions of a (primitive) set-optimization problem, as introduced in \cite{HeydeLoehne11} and deepened in \cite{Loehne11Book}. To prove the main results we need also to deal with scalarization problems. However, while in the vector case this might only be a technical need, we prove that eventually the set-valued variational inequalities and their scalar counterparts provide different insights on the problem. Some relevant information on the solution of the set-optimization problem is provided only through the scalar version of the inequality. The special case of vector optimization is finally studied, to recover classical results stated in \cite{CreGinRoc08, YangYang04}.\\
The paper is organized as follows. Section 2 is devoted to preliminary results on set-optimization that will be used throughout the paper. The concept of solution to a set-optimization problem and the Dini-type derivatives are presented and some properties are proved. Section 3 presents the main results. As the solution concept relies on two properties, we develop two different sets of relations between our variational inequalities and the set-optimization. The first one provides a variational characterization of 'infimizer', while the second one is devoted to characterize 'minimizer'.
Finally, Section 4 applies the previous results to vector optimization. 
The relations proved for the convex case in this paper reproduce those already known for the vector case between optimization and variational inequalities. We leave as an open question whether convexity can be relaxed, as indeed can be done for vector-valued functions.


\section{Preliminaries}

\subsection{Order and operations with sets}

Throughout the paper, unless explicitly stated otherwise, we assume the setting and notation introduced in this section.

Let $Z$ be a locally convex Hausdorff space with  dual space $Z^*$. The set $\mathcal U$ is the set of all closed, convex and balanced $0$-neighborhoods in $Z$, a $0$-neighborhood base of $Z$. By $\cl A$, $\co A$ and $\Int A$, we denote the closed or convex hull of a set $A\subseteq Z$ and the topological interior of $A$, respectively. The conical hull of a set $A$ is $\cone A=\cb{ta\st a\in A,\, 0<t}$.

The recession cone of a nonempty closed convex set $A\subseteq Z$ is given by
\begin{equation}\label{eq:Rec_cone}
0^+A=\cb{z\in Z\st A+\cb{z}\subseteq A},
\end{equation}
a closed convex cone, \cite[p.6]{Zalinescu02}. By definition, $0^+\emptyset=\emptyset$ is assumed.

$Z$ is pre-ordered through a closed convex cone $C$ with $0\in C$ and nontrivial negative dual cone
\[
C^-=\cb{z^*\in Z^*\st \forall c\in C:\; z^*(c)\leq 0},
\] 
$\C\neq \emptyset$ by setting
\[
z_1\leq z_2\quad\Leftrightarrow\quad \cb{z_2}+C\subseteq \cb{z_1}+C
\]
for all $z_1, z_2\in Z$.  This relation is extended to $\P(Z)$, the power set of $Z$ including $\emptyset$ and $Z$ (compare \cite{Hamel09}
and the references therein) by setting
\[
A_1\lel A_2\quad\Leftrightarrow\quad A_2+C\subseteq A_1+C
\]
for all $A_1, A_2\subseteq Z$, . 

We introduce the subset
\[
\G(Z,C)=\cb{A\subseteq Z\st A=\cl\co(A+C)}
\]
which is an order complete lattice and $A_1\lel A_2$ is equivalent to $A_1\supseteq A_2$ whenever $A_1, A_2\in\G(Z,C)$ .
For any subset $\mathcal A\subseteq \G(Z,C)$, supremum and infimum of $\mathcal A$ in $\G(Z,C)$ are given by
\begin{align}\label{eq:inf_sup}
\inf\mathcal A=\cl\co \bigcup\limits_{A\in \mathcal A}A;\quad \sup\mathcal A=\bigcap\limits_{A\in \mathcal A}A
\end{align}
and for a net $\cb{A_i}_{i\in I}$ in $\G(Z,C)$, limit inferior and limit superior are defined accordingly,
\begin{align}\label{eq:liminf_sup}
\liminf A_i=\bigcap\limits_{j\in I}\cl\co \bigcup\limits_{i\geq j}A_i;\quad \limsup\mathcal A=\cl\co\bigcup\limits_{j\in I}\bigcap\limits_{i\geq j}A_i.
\end{align}

When $\mathcal A=\emptyset$, then we agree on $\inf\mathcal A=\emptyset$ and $\sup\mathcal A =Z$. Especially, $\G(Z,C)$ possesses a greatest and smallest element $\inf\G(Z,C)=Z$ and $\sup\G(Z,C)=\emptyset$.

%

The Minkowsky addition and multiplication with negative reals need to be slightly adjusted to provide operations on $\G(Z,C)$. We define
\begin{align}
\forall A, B\in \G(Z,C):\quad 		&A\oplus B=\cl\cb{a+b\in Z\st a \in A,\, b\in B};\\
\forall A\in\G(Z,C),\, \forall 0<t:\quad &t\cdot A=\cb{ta \in Z\st a\in A };\\
\forall A\in\G(Z,C):\quad 			&0\cdot A=C.
\end{align}
Especially, $0\cdot \emptyset=0\cdot Z=C$ and $\emptyset$ dominates the addition in the sense that $A\oplus \emptyset=\emptyset$ is true for all $A\in\G(Z,C)$.
Moreover, $A\oplus C=A$ is satisfied for all $A\in \G(Z,C)$, thus $C$ is the neutral element with respect to addition.

As a consequence, 
\begin{align}
\forall \mathcal A\subseteq \G(Z,C),\, \forall B\in\G(Z,C) :\quad B\oplus \inf\mathcal A=\inf\cb{B\oplus A\st A\in \mathcal A},
\end{align}
or, equivalently,  the $\inf$--residual
\begin{align}
A\idif B=\inf\cb{M\in\G(Z,C)\st A\lel B\oplus M}
\end{align}
exists for all $A, B\in \G(Z,C)$. The following properties are well known in lattice theory, compare also \cite[Theorem 2.1]{HamelSchrage12}.
\begin{align}
A\idif B    &=\cb{z\in Z\st B+\cb{z}\subseteq A};\\
A             &\lel B\oplus (A\idif B)
\end{align}

Overall, the structure of $\G^\triup=\of{\G(Z,C),\oplus,\cdot,C,\lel}$ is that of an $\inf$--residuated conlinear space, compare also   \cite{Fuchs66}, \cite{galatos2007residuated}, \cite{GetanMaLeSi}, \cite{HamelHabil} \cite{MartinezLegazSinger95}.


Historically, it is interesting to note that R. Dedekind \cite{Dedekind1872} introduced the residuation concept and used it in order to construct the real numbers. The construction above is in this line of ideas, but in a rather abstract setting.

\begin{example}
\label{ExExtReals}
Let us consider $Z = \R$, $C = \R_+$. Then $\G\of{Z, C} = \cb{[r, +\infty) \mid r \in \R}\cup\cb{\R}\cup\cb{\emptyset}$, and $\G^\triup$ can be identified (with respect to the algebraic and order structures which turn $\G\of{\R, \R_+}$ into an ordered conlinear space and a  complete lattice admitting an inf-residuation) with $\OLR = \R\cup\cb{\pm\infty}$ using the 'inf-addition' $\isum$ (see \cite{HamelSchrage12}, \cite{RockafellarWets98}). The inf-residuation on $\OLR$ is given by
\[
r \idif s  = \inf\cb{t\in\R \mid r \leq s \isum t}
\]
for all $r,s\in\OLR$, compare \cite{HamelSchrage12} for further details.
\end{example}

Each element of $\G^\triup$ is closed and convex and $A=A+C$. Hence, by a separation argument we can prove
\begin{align}\label{eq:scal_representation_Set}
\forall A\in \G^\triup:\quad A=\bigcap\limits_{z^*\in\C}\cb{z\in Z\st -\sigma(z^*| A)\leq -z^*(z)},
\end{align}
where $\sigma(z^*| A)=\sup\cb{z^*(z)\st z\in A}$ is the support function of $A$ at $z^*$. Especially, $A=\emptyset$ if and only if there exists $z^*\in \C$ such that $-\sigma(z^*| A)=+\infty$.

\begin{lemma}\label{prop:scal_of_infimum}\cite[Proposition 3.5]{Schrage10Opt}
Let $\mathcal A\subseteq\G^\triup$ be a given subset, then
\begin{align}
&\inf\mathcal A=\bigcap\limits_{z^*\in\C}\cb{z\in Z\st \inf\cb{-\sigma(z^*| A)\st A\in\mathcal A}\leq -z^*(z)}\\
&\forall z^*\in \C:\quad -\sigma(z^*| \inf\mathcal A)=\inf\cb{-\sigma(z^*| A)\st A\in\mathcal A}.
 \end{align}
\end{lemma}

\begin{lemma}\label{lem:scal_of_difference}\cite[Proposition 5.20]{HamelSchrage12}
Let $A, B\in\G^\triup$, then
\begin{align}
&A\idif B=\bigcap\limits_{z^*\in \C}\cb{z\in Z\st (-\sigma(z^*| A))\idif (-\sigma(z^*| B))\leq -z^*(z)};\\
&\forall z^*\in \C:\quad \of{-\sigma(z^*| A)}\idif\of{-\sigma\of{z^*|B}}\leq -\sigma(z^*| A\idif B).
\end{align}
\end{lemma}

In general, the difference of the scalarizations and the scalarization of the difference do not coincide, as the following example shows.
\begin{example}\label{ex:difference_can_be_empty}
Let $Z=\R^2$ and $C=\cl\cone{(0,1)^T}$,  $B=\cb{(x,y)\in \R^2\st -1\leq x\leq 1,\, 0\leq y}$ and $A=C$. Then $(-\sigma(z^*| A))\idif (-\sigma(z^*| B))\in\R$ is satisfied for all $z^*\in \C$ and
\begin{align*}
A\idif B=\cb{z\in Z\st 1\leq (-1,0)^T z,\, 1\leq (1,0)^T z,\, 0\leq (0,1)^T z }=\emptyset,
\end{align*}
 thus $-\sigma(A\idif B)=+\infty$.
\end{example}

The following rules  will be used frequently later on.
\begin{lemma}\label{lem:calc_of_diff}
Let $A,B,D\in\G^\triup$, $0<s$ and $t\in\of{0,1}$ be given, then
\begin{enumerate}[(a)]
\item
$s(A\idif B)=sA\idif sB;$
\item
$(tA\oplus (1-t)B)\idif D\lel t(A\idif D)\oplus (1-t)(B\idif D);$
\item
$A\idif D\lel\of{A\idif B}\oplus\of{B\idif D};$
\item
If $A\neq \emptyset$, then $0^+A=\of{A\idif A}$.
\end{enumerate}
\end{lemma}
\proof
\begin{enumerate}[(a)]
\item 
It holds $z\in (A\idif B)$ if and only if $B+\cb{z}\subseteq A$ or equivalently 
$sA\lel sB+\cb{sz}$.
\item 
As $D\in \G^\triup$ is assumed, $tD\oplus (1-t)D=D$. Let $z_A\in A\idif D$ and $z_B\in B\idif D$ be given, then $ (tA\oplus (1-t)B)\lel D+(tz_A+(1-t)z_B)$ is satisfied.
\item 
The inclusion is true if and only if
\[
A\lel \of{A\idif B}\oplus\of{B\idif D}\oplus D.
\]
As we know that $B\lel \of{B\idif D}\oplus D$ and $A\lel \of{A\idif B}\oplus B$, this inclusion is true.
\item This is immediate from the definition of $0^+A$.
\end{enumerate}
\pend

Lemma \ref{lem:calc_of_diff} (d) suggests that, if needed, we can use the recession cone of a set as $0$--element in certain inequalities. It is remarkable that   for any $A\in\G^\triup$, either $A=\emptyset$, or $ 0^+A\lel C$. 
To implement these remarks in the sequel, we use the following properties of recession cones.

\begin{proposition}\label{prop:Rec_A}
Let $A\in\G^\triup\setminus\cb{\emptyset}$, then
\begin{equation*}
0^+A=\cb{z\in Z\st \forall z^*\in \C:\; -\sigma(z^*|A)=-\infty\,\vee\, 0\leq -z^*(z)}.
\end{equation*}
Especially, for all $A\in\G^\triup$, either $A=\emptyset$, or
\begin{equation}\label{eq:RecCone_finiteScal}
0^+ A=\bigcap\limits_{\substack{z^*\in\C\\ -\sigma(z^*| A)\in\R}}\cb{z\in Z\st 0\leq -z^*(z)}.
\end{equation}
\end{proposition}
\proof
Assume $z\notin 0^+A$, then either $A=\emptyset$ or there exists a $z^*\in Z^*$ such that $\sigma(z^*|A)<z^*(a+z)$ is satisfied for some $a\in A$. As $z^*(a+z)\leq \sigma(z^*|A)+z^*(z)$, this implies   
$-z^*(z)< 0$ and $-\sigma(z^*|A)\neq-\infty$ and therefore $z^*\in\C$.
On the other hand, assume $z\in 0^+A$, then $A$ is nonempty and $A+\cb{z}\subseteq A$, hence for all $z^*\in Z^*$ it holds
$\sigma(z^*|A+\cb{z})\leq \sigma(z^*|A)$, hence $\sigma(z^*|A)+z^*(z)\leq \sigma(z^*|A)$. This implies that either $-\sigma(z^*|A)=-\infty$ or $0\leq-z^*(z)$ is true for all $z^*\in Z^*$ and thus especially for $z^*\in\C$.

If $A=Z$, then $-\sigma(z^*|Z)=-\infty\notin\R$ is satisfied for all $z^*\in\C$, hence \eqref{eq:RecCone_finiteScal} is true with $0^+Z=Z$.
Hence let $A\neq Z$ or $\emptyset$, then $-\sigma(z^*| A)\notin\R$ implies $-\sigma(z^*| A)=-\infty$ and the statement is proved.

\pend

\begin{lemma}
Let $A\in\G^\triup\setminus\cb{\emptyset}$, then
\[
\cb{z^*\in Z^*\st -\sigma(z^*|A)\in\R}
\subseteq (0^+A)^-\subseteq C^-.
\]
\end{lemma}
\proof
Assume $-\sigma(z^*|A)\in\R$ and $A+\cb{z}\subseteq A$. Then 
\[
-\sigma(z^*|A)\leq -\sigma(z^*|A+\cb{z})=-\sigma(z^*|A)+(-z^*(z))
\]
implies $0\leq -z^*(z)$, in other words $z^*\in (0^+A)^-$.
The second inclusion is immediate, as $A\in\G^\triup\setminus\cb{\emptyset}$ implies $0^+A\supseteq C$.
\pend

\begin{lemma}\label{lem:rec_cone_monotone}
Let $A, B\in \G^\triup\setminus\cb{\emptyset}$, then
\begin{eqnarray*}
&0^+(A\oplus B)\lel \cl\co\of{0^+A\cup 0^+B}= 0^+A \oplus 0^+B;\\
&A\lel B\quad\Rightarrow\quad  0^+A\lel 0^+B.
\end{eqnarray*}
\end{lemma}
\proof
Assume $A+\cb{z_A}\subseteq A$ and $B+\cb{z_B}\subseteq B$, then for all $a\in A$ and all $b\in B$ it holds
\[
a+b+(z_A+z_B)\in A\oplus B
\]
and as both $0^+A$ and $0^+B$ are convex cones, for all $t\in\sqb{0,1}$ it holds
\[
ta+(1-t)b+(z_A+z_B)\in A\oplus B.
\]
If $z\in A\oplus B$, then for all $U\in\mathcal U$ there exist $a\in A$, $b\in B$ and $t\in\sqb{0,1}$ with $ta+(1-t)b\in \cb{z}+U$, such that
\[
ta+(1-t)b+(z_A+z_B)\in \cb{z+(z_A+z_B)}+U,
\] 
and hence $z+(z_A+z_B)\in A\oplus B$, proving $0^+A+0^+B\subseteq 0^+(A\oplus B)$.
As $A\oplus B$ is a closed convex set, the recession cone is a closed convex cone, so
\[
0^+A\oplus 0^+B=\cl\co(0^+A+0^+B)\subseteq 0^+(A\oplus B).
\]
Since $0\in 0^+A\cap 0^+B$ implies $0^+A\cup 0^+B\subseteq 0^+A \oplus 0^+B$, also  $\cl\co\of{0^+A\cup 0^+B}\subseteq 0^+A \oplus 0^+B$ holds true.
On the other hand, if $z_A\in 0^+A$ and $z_B\in 0^+B$ are given, then $z_A+z_B\in\co\of{0^+A\cup 0^+B}$, hence
$\cl\co\of{0^+A\cup 0^+B}\supseteq 0^+A \oplus 0^+B$ proves equality

Finally, let $A\lel B$ be satisfied, $B+\cb{z}\subseteq B$ and $a+z\notin A$ for some $a\in A$.
Then there exists a neighborhood $U\in\mathcal U$ such that $\cb{a+z}+U\cap A=\emptyset$, as $A$ is closed and thus there is exists $t\in\of{0,1}$ such that
\[
t\of{b+\frac{1}{t}z}+(1-t)a=a+z+t(b-a)\in \cb{a+z}+U.
\]
But since $A$ is convex and $0^+B$ is a cone, this implies 
\[
t\of{b+\frac{1}{t}z}+(1-t)a\in \co(B+A)\subseteq A,
\]
a contradiction.
\pend

Moreover, we can remark that  for any $A\in \G^\triup$ the following properties hold true
\begin{enumerate}[(i)]
\item $0^+A \oplus 0^+\emptyset= 0^+(A\oplus \emptyset)$; 
\item $0^+A\lel 0^+\emptyset$.
\end{enumerate}
On the contrary, 
$ 0^+A \oplus 0^+\emptyset\lel 0^+A\cup 0^+\emptyset$
can be proven if and only if $A=\emptyset$.

\begin{lemma}\label{lem:Rec_A2}
If $A\idif B\neq \emptyset$, then
\begin{equation*}
0^+(A\idif B)\lel 0^+A\lel 0^+B.
\end{equation*}
If additionally $B\neq \emptyset$, then we also get
\begin{equation*}
0^+(A\idif B)= 0^+A.
\end{equation*}
\end{lemma}
\proof
Assume $A\idif B\neq \emptyset$.
If $B=\emptyset$, then $A\idif B=Z$ and the first equation is immediate. Hence let $B\neq\emptyset$.
Then
$\emptyset\neq B\oplus (A\idif B)\subseteq A$ and because $A$ is closed and convex by assumption, we can apply Lemma \ref{lem:rec_cone_monotone} to prove
\[
0^+B \cup 0^+(A\idif B)\subseteq 0^+(B\oplus\of{A\idif B})\subseteq 0^+A.
\] 
On the other hand, if $B+\cb{z}\subseteq A$, that is $z\in A\idif B$, then for all $z_0\in 0^+A$ it holds
$B+\cb{z+z_0}\subseteq A$, hence $0^+A\subseteq 0^+(A\idif B)$.
\pend

%
%
%




\subsection{Set-valued functions}

Let $X$ be a linear space. A function $f \colon X \to
\G^\triup$ is called convex when
\begin{equation}
\label{EqConvFct} \forall x_1, x_2 \in X, \; \forall t \in \of{0,1} \colon
 f\of{tx_1+(1-t)x_2} \lel tf\of{x_1} \oplus \of{1-t}f\of{x_2}.
\end{equation}
It is an easy exercise (see, for instance, \cite{Hamel09}) to show that $f$ is convex if and only if the set
\[
\gr f = \cb{\of{x,z} \in X \times Z \colon z \in f\of{x}}
\]
is convex. A $\G^\triup$-valued function $f$ is called positively homogeneous when
\[
\forall 0<t, \forall x \in X \colon f\of{tx} \lel tf\of{x},
\]
and it is called sublinear if it is positively homogeneous and convex. It can be shown that $f$ is sublinear if and only if $\gr f$ is a convex cone. Compare also \cite[Definition 2.1.1.]{AubinFrankowska} on above definitions.

The (effective) domain of a function $f:X\to\G^\triup$ is the set
$\dom f=\cb{x\in X\st f(x)\neq \emptyset}$.
Since $\emptyset$ is the supremum of $\G^\triup$, the previous notion of domain of a set-valued function extends the scalar notion of effective domain.
The image set of a subset $A\subseteq X$ through $f$  is denoted by 
\[
f\sqb{A}=\cb{f(x)\in\G^\triup\st x\in A}\subseteq \G^\triup.
\]
We underline that $f\sqb{A}$ is a subset of $\P(Z)$ rather then a subset of $Z$, while $\inf f\sqb{A}=\cl\co\bigcup\limits_{a\in A}f(a)$ is an element of $\P(Z)$, hence a subset of $Z$. 

\begin{proposition}\label{prop:Rec_A3}
Let $f:X\to\G^\triup$ be convex, $x_0\in\dom f$. If $x\in\dom f$, then  $t\mapsto 0^+(f(x+t(x_0-x)))$ is constant on $\of{0,1}$ and
$0^+(f(x+t(x_0-x)))\lel 0^+f(x)\cup 0^+f(x_0)$ is satisfied  for all $t\in\of{0,1}$.
\end{proposition}
\proof
Let $t\in\sqb{0,1}$ and denote $x_t=x+t(x_0-x)$. By convexity of $f$, for any $z_0\in 0^+f(x_0)$ and $z\in 0^+f(x)$, $z_t=tz+(1-t)z_0\in 0^+f(x_t)$ is satisfied. Since both recession cones contain $0$, we have $z_0+0\in 0^+f(x_t)$ and $z+0\in 0^+f(x_t)$. Therefore $0^+f(x_t)\supseteq 0^+f(x_0)\cup 0^+f(x)$. 

Moreover let $0<s<t<1$ be given. By replacing $x$ with $x_t$ in above argument we have 
\[
0^+f(x_s)\supseteq 0^+f(x_0)\cup 0^+f(x_t)=0^+ f(x_t)
\]
 and by replacing $x_0$ by $x_s$ instead we have 
\[
0^+f(x_t)\supseteq 0^+f(x_s)\cup 0^+f(x)=0^+ f(x_s),
\]
hence $0^+ f(x_s)=0^+f(x_t)$ is proven for all $s,t\in\of{0,1}$.
\pend

Given a function $f \colon X \to \G^\triup$, the family of extended real-valued functions $\vp_{f, z^*} \colon  X \to \R\cup\cb{\pm\infty}$ defined by
\[
\vp_{f,z^*}\of{x} = \inf\cb{-z^*\of{z} \mid z \in f\of{x}}, \; z^* \in C^-\bs\cb{0}
\]
is the family of scalarizations of $f$. 
Some properties of $f$ are inherited by its scalarizations and vice versa. For instance, 
$f$ is convex if and only if  $\vp_{f, z^*}$ is convex for each $z^* \in C^-\bs\cb{0}$. In turn, convexity of $\vp_{f,z^*}$ is equivalent to convexity of the function $f_{z^*}:X\to\G^\triup$ given by 
\[
f_{z^*}(x) = \cb{z \in Z \st \vp_{f, z^*}\of{x} \leq -z^*\of{z}}.
\]
Moreover, a standard separation argument shows
\[
\forall x \in X \colon f\of{x} = \bigcap_{z^* \in C^-\bs\cb{0}}f_{z^*}(x).
\]

\begin{remark}
The function $f_{z^*}:X\to\G^\triup$ maps $x$ to the sublevel set $L^\leq_{z^*}(-\vp_{f,z^*}(x))$ of $z^*$ to the level $-\vp_{f,z^*}(x)$. For all $z^*\in\C$ and all $x\in X$ it holds
\begin{equation}
f_{z^*}(x)=L^\leq_{z^*}(\sigma(z^*|f(x)))=\cb{z\in Z\st z^*(z)\leq -\vp_{f,z^*}(x)}.
\end{equation}
Therefore either $f_{z^*}(x)\in\cb{\emptyset, Z}$, or it is a closed affine half space with a supporting point $z\in f_{z^*}(x)$ such that $\vp_{f,z^*}(x)=-z^*(z)$.
If $f(x)\neq \emptyset$, then either $f_{z^*}(x)=Z$, or $\vp_{f,z^*}(x)\in\R$, thus
\begin{align*}
\forall x\in X:\quad f(x)=\emptyset\;\vee\; f(x)=\bigcap\limits_{\substack{z^*\in\C: \\ \vp_{f,z^*}(x)\in\R}}f_{z^*}(x).
\end{align*}
\end{remark}

%

\begin{definition}\label{def:l.s.c.}
\begin{enumerate}[(a)]
\item 
Let $\vp:X\to\OLR$ be a function, $x_0\in X$. Then $\vp$ is said to be lower semicontinuous (l.s.c.) at $x_0$, iff  
\begin{equation*}
\forall r\in\R:\quad r<\vp(x_0)\;\Rightarrow\;
\exists U\in\mathcal U:\, \forall u\in U:\,  r<\vp(x_0+u).
\end{equation*} 

\item 
Let $f:X\to\G^\triup$ be a function, $M^*\subseteq \C$. Then $f$ is said $M^*$-- lower semicontinuous ($M^*$--l.s.c.) at $x_0$, iff
$\vp_{f,z^*}$ is l.s.c. at $x_0$ for all $z^*\in M^*$.
\item 
Let $f:X\to\G^\triup$ be a function. If 
\[
f(x)\lel \liminf\limits_{u\to 0}f(x+u)=\bigcap\limits_{U\in\mathcal U}\cl\co \bigcup\limits_{u\in U}f(x+u)
\]
is satisfied, then $f$ is lattice lower semicontinuous (lattice l.s.c.) at $x$.
\item
A function $f:X\to\G^\triup$ is lattice l.s.c. iff it it is lattice l.s.c. everywhere.
\end{enumerate}
\end{definition}


In \cite{HeydeSchrage11R}, it has been proven that if $f$ is $\of{\C}$--l.s.c. at $x$, then it is also lattice l.s.c. at $x$.
One can show that if $f$ is convex, then $f$ is lattice l.s.c. if and only if $\gr f=\cb{(x,z)\st z\in f(x)} \subseteq X \times Z$ is a closed set with respect to the product topology, see \cite{HamelSchrage13PJO}. 

In \cite{HeydeSchrage11R}, a detailed study of continuity concepts for set valued functions is proposed. Indeed it is also shown that none of the concepts in Definition \ref{def:l.s.c.} coincides with those used in some literature (see e.g.  \cite{AliprantisBorder, AubinFrankowska, GRTZ}).

\begin{remark}
For notational simplicity the restriction of a set-valued function $f:X\to\G^\triup$ to a segment with end points $x_0,x\in X$ is denoted by $f_{x_0,x}:\R\to\G^\triup$ with
\[
f_{x_0,x}(t)=\begin{cases}
	f(x_0+t(x-x_0)),\text{ if } t\in\sqb{0,1};\\
	\emptyset,\text{ elsewhere.}
\end{cases}
\]
Equivalently, the restriction of a scalar-valued function $\vp:X\to\OLR$ to the same segment is defined by
\[
\vp_{x_0,x}(t)=\begin{cases}
	\vp(x_t),\text{ if } t\in\sqb{0,1};\\
	+\infty,\text{ elsewhere.}
\end{cases}
\]
Setting $x_t=x_0+t(x-x_0)$ for all $t\in\R$,
the scalarization of the restricted function $f_{x_0,x}$ is equal to the restriction of the scalarization of $f$ for all $z^*\in\C$.

If $f$ is convex, $x_0,x_t\in \dom f$ for some $t\in\of{0,1}$, then $\of{\vp_{f,z^*}}_{x_0,x}$ is lower semicontinuous on $\of{0,t}$, hence $f_{x_0,x}$ is lattice l.s.c. on $\of{0,t}$.
\end{remark}


The following notion, introduced in \cite{HamelSchrage13PJO}, is used in the sequel.

\begin{definition}\label{def:CanExt}
Let $f \colon X \to \G^\triup$ be a function and $M \subseteq X$. We define the
inf-translation of $f$ by $M$ to be the function $\hat f\of{\cdot; M} \colon X \to
\G^\triup$ given by
\begin{equation}\label{EqSetTranslation}
\hat f\of{x; M} = \inf f\sqb{M +\cb{x}} = \cl\co\bigcup_{m \in M}f\of{m+x}.
\end{equation}
\end{definition}

The function $\hat f\of{\cdot; M}$ is nothing but the canonical extension of $f$ at
$M + \cb{x}$ as defined in \cite{HeydeLoehne11}. 
The following properties of the inf-translation are used in the proofs of the main results.  

\begin{lemma}\label{lem:fcoM_convex}\cite[Lemma 5.8 (b)]{HamelSchrage13PJO}
Let $f:X\to\G^\triup$ be convex, $M\subseteq X$, then
$\hat f\of{\cdot; \co M}:X\to\G^\triup$ is convex.
\end{lemma}

\begin{lemma}\label{lem:scal_of_fM}
Let $f:X\to\G^\triup$, $z^*\in\C$ and $M\subseteq X$ be nonempty. Then
\begin{equation*}
\forall x\in X:\quad \inf \vp_{f,z^*}\sqb{M+\cb{x}}=\vp_{\hat f\of{\cdot; M},z^*}(x).
\end{equation*} 
Moreover, by defining $\hat\vp_{f,z^*}(x;M)=\inf \vp_{f,z^*}\sqb{M+\cb{x}}$, it holds
\begin{equation*}
\forall x\in X:\quad \hat\vp_{f,z^*}(x;M)=\vp_{\hat f\of{\cdot; M},z^*}(x),
\end{equation*} 
that is the operations of taking the inf translation of a function and taking its scalarization commute.
\end{lemma}
\proof
The statement is an easy consequence of Lemma \ref{prop:scal_of_infimum}.
\pend

\begin{lemma}\label{lem:dom_of_fM}
Let $f:X\to\G^\triup$ and $M\subseteq X$ be nonempty, then the domain of $\hat f\of{\cdot; M}:X\to\G^\triup$ is the set
\begin{equation}
\dom \hat f\of{\cdot; M}=\bigcup\limits_{m\in M}\dom f+\cb{-m}.
\end{equation}
\end{lemma}
\proof
Since $x\in\dom \hat f(\cdot;M)$ if and only if $\inf f\sqb{M+\cb{x}}\neq \emptyset$, there exists $m\in M$ such that $f(m+x)\neq\emptyset$.
Therefore, $x\in\dom \hat f(\cdot;M)$ if and only if  $m+x\in\dom f$ for some $m\in M$. In other words
$x\in \bigcup\limits_{m\in M}\dom f+\cb{-m}$.
\pend

\begin{lemma}\label{lem:f_M_is convex,lsc}
Let $f:X\to\G^\triup$ be convex, $M\subseteq X$ a nonempty set and $z^*\in\C$, If any of the following conditions is satisfied, then the restriction of  $\hat f(\cdot;\co M)$ to the segment $\sqb{0,x}$ is $\of{\C}$--l.s.c. at $0$ for all $x\in X$.
\begin{enumerate}[(a)]
\item
$\hat f(0;M)=\inf f\sqb{X}$;
\item $0\in \Int \bigcup\limits_{m\in \co M}(\dom f+\cb{-m})$;
\item  $\of{\vp_{f,z^*}}_{m,x}:X\to\OLR$ is continuous at $0$ for all $m\in\co M$, $x\in X$ and all $z^*\in \C$.
\end{enumerate}
\end{lemma}
\proof
\begin{enumerate}[(a)]
\item
If $\hat f(0;M)=\inf f\sqb{X}$, then 
$\vp_{\hat f(\cdot;\co M),z^*}(0)=\inf \vp_{\hat f(\cdot;\co M),z^*}\sqb{X}$ is true for all $z^*\in\C$. Hence each scalarization $\vp_{\hat f(\cdot;\co M),z^*}$ is l.s.c. at $0$ and therefore  $\hat f(\cdot;\co M)$ is $\C$--l.s.c at $0$.
\item
By Lemma \ref{lem:dom_of_fM}, $\bigcup\limits_{m\in \co M}(\dom f+\cb{-m})$ is the domain of $\hat f(\cdot;\co M)$ and by Lemma \ref{lem:fcoM_convex}, $\hat f(\cdot;\co M)$ is convex. This is true if and only if each scalarization of $\hat f(\cdot;\co M)$ i.e. $\of{\hat\vp_{f,z^*}}(\cdot;\co M)$ is convex, compare Lemma \ref{lem:scal_of_fM}. If $0\in \Int \bigcup\limits_{m\in \co M}(\dom f+\cb{-m})$ is assumed, then the restriction of each scalarization $\vp_{f,z^*}(\cdot;\co M)$ to $\sqb{x_0,x}$ is l.s.c. at $0$,  as $\dom \hat f(\cdot;\co M)=\dom \of{\hat\vp_{f,z^*}}\of{\cdot;\co M}$. 
\item 
Let $\of{\vp_{f,z^*}}_{m,x}:X\to\OLR$ be continuous at $0$ for all $m\in\co M$ and all $x\in X$.
In this case,
\begin{align*}
\limsup\limits_{t\downarrow 0} (\vp_{\hat f(\cdot;\co M),z^*})_{0,x}(t)
		&=\limsup\limits_{t\downarrow 0}\inf\limits_{m\in\co M}\of{\vp_{f,z^*}}_{m,x}(t)\\
		&\leq \inf\limits_{m\in\co M}\limsup\limits_{t\downarrow 0}\of{\vp_{f,z^*}}_{m,x}(t)\\
		&= \inf\limits_{m\in\co M}\of{\vp_{f,z^*}}_{m,x}(0)\\
		&=\hat\vp_{f,z^*}\of{0;\co M}.
\end{align*}
Hence for each $z^*\in\C$, the restriction of  $\vp_{f,z^*}(\cdot;\co M)$ to $\sqb{0,x}$ is convex and u.s.c. at $0$, therefore l.s.c. at $0$, too.
\pend
\end{enumerate}



In this framework, 
we are interested to study the problem
\begin{align}\label{eq:optim_problem}
\tag{P} \mbox{minimize} \quad f(x) \quad \mbox{subject to} \quad  x \in X
\end{align}
where $f$ is a $\G^\triup$-valued function. Following \cite{HeydeLoehne11}, to solve \eqref{eq:optim_problem} means to  look for both
the infimum in $\G^\triup$, as introduced in \eqref{eq:inf_sup}, and for subsets of $X$ where the infimum is attained. This approach is different from most other approaches in set optimization, see for example \cite[Definition 14.2]{Jahn04}, \cite{HernandezRodriguezMarin07-1}, \cite{HernandezRodriguezMarin07-2} and the references therein.

More formally, the solution concept based on Definitions \ref{DefInfimizer} and \ref{def:Minimizer} is stated in Definition \ref{def:Infimizer_Minimizer_Solution}.

\begin{definition}
\label{DefInfimizer}
 Let $f \colon X \to \G^\triup$. A subset $M \subseteq X$
is called an {\em infimizer} of $f$ if
\[
\inf \cb{f(m)\st m\in M}=\inf \cb{f(x)\st x\in X}.
\]
\end{definition}

According to the definition of $\hat f\of{\cdot; M}:X\to\G^\triup$, it follows easily that
\begin{equation*}
\forall M\neq\emptyset:\quad \inf \cb{\hat f\of{x; M}\st x\in X}=\inf \cb{f(x)\st x\in X}
\end{equation*}
and $M$ is an infimizer of $f$ if and only if $\cb{0}$ is an infimizer of $\hat f\of{\cdot; M}:X\to\G^\triup$,
\begin{equation*}
\hat f\of{0; M}=\inf \cb{\hat f\of{x; M}\st x\in X} \quad\Leftrightarrow\quad  \inf \cb{f(m)\st m\in M}=\inf \cb{f(x)\st x\in X}.
\end{equation*}

\begin{proposition}\cite[Proposition 5.9]{HamelSchrage13PJO}\label{prop:infimizer}
Let $f:X\to\G^\triup$ be convex and $M\subseteq X$, then the following are equivalent
\begin{enumerate}[(a)]
\item $M$ is an infimizer of $f$;
\item $\cb{0}$ is an infimizer of $\hat f(\cdot;M)$;
\item  $\cb{0}$ is an infimizer of $\hat f(\cdot;\co M)$ and $\hat f(0;M)=\hat f(0;\co M)$. 
\end{enumerate}
\end{proposition}



\begin{proposition}\label{prop:infimizer_iff}
Let $f:X\to\G^\triup$ and $x_0\in\dom f$. Then the following are equivalent
\begin{enumerate}[(a)]
\item 
$f(x_0)=\inf f\sqb{X}$;
\item 
$
\forall x\in X,\,\forall z^*\in\C:\,  \vp_{f,z^*}(x_0)\leq \vp_{f,z^*}(x); 
$
\item 
$
\forall x\in X,\,\forall z^*\in\C:\,  
\vp_{f,z^*}(x_0)\idif \vp_{f,z^*}(x)\leq 0; 
$
\item
$
\forall x\in X:\, 0\in f(x_0)\idif f(x).
$
\item 
$
\forall x\in X,\, \forall z^*\in\C:\,  
\vp_{f,z^*}(x_0)=-\infty\,\vee\,0\leq \vp_{f,z^*}(x)\idif \vp_{f,z^*}(x_0). 
$
\end{enumerate}
Each of these conditions implies
\begin{enumerate}[(f)]
\item
$
\forall x\in X:\, 0^+f(x_0)\lel f(x)\idif f(x_0).
$
\end{enumerate}
\end{proposition}
\proof
The equivalence between (a), (b), (c) and (e) is immediate. By Lemma \ref{lem:scal_of_difference} (c) and (d) are equivalent and by Proposition \ref{prop:Rec_A}, (e) implies (f).
\pend

\begin{remark}
For scalars $a,b\in\R$, $a\le b$ can be equivalently stated as $a-b\le 0$ or $0\le b-a$. For $A,\,B\in \G^\triup\setminus\left\{\emptyset\right\}$ we have a similar result for the equivalence between $A\lel B$ and $A\idif B\lel C$ (and actually as '$A\idif B\lel  0^+A$' or $0\in A\idif B$).\\
On the other hand, $A\lel B$ only implies $0^+A\lel B-A$. Moreover $0^+B$ is not necessarily equal to $C$, the neutral element in $\G^\triup$, but $0^+A\lel C$, whenever $A\neq\emptyset$.


\end{remark}
As $\dom f$ is always an infimizer of $f$, further requirements are usually assumed on the values $f\of{x}$, $x \in M$, for $M$ to be a solution.  e.g. $f(x)$ is minimal in some sense, compare
\cite{HamelSchrage13PJO, HeydeLoehne11, Loehne11Book}.

\begin{definition}\label{def:Minimizer}
Let $f:X\to \G^\triup$ be given. An element $x_0\in X$ is called a {\em minimizer} of $f$, iff $f(x_0)$ is minimal in $f\sqb{X}$, i.e.
\begin{align}\label{eq:Min}
\tag{Min}
\forall x\in X:\quad f(x)\lel f(x_0)\quad\Rightarrow\quad f(x)=f(x_0).
\end{align}
The set of all minimal elements of $f\sqb{X}$ is denoted by $\Min f\sqb{X}$.
\end{definition} 
If $x_0$ is a minimizer of a convex (set-valued) function $f$, then $f(x)=f(x_0)$ is satisfied if and only if $f$ is constant on the set $\cb{x_t\in X\st x_t=x_0+t(x-x_0),\, t\in\sqb{0,1}}$.

Notice that if $M=\cb{x}$ is an infimizer, then $x$ automatically is a minimizer of $f$. On the other hand, a set of minimizers is not necessarily an infimizer. Let  $\psi:S\subseteq X\to Z$ and its epigraphical extension $f=\psi^C:X\to \G^\triup$, defined by 
\begin{align}\label{eq:sv_extension}
f(x)=\begin{cases}
	\cb{\psi(x)}+C, &\text{ if } x\in S;\\
	\emptyset, &\text{ elsewhere.}
	\end{cases}
\end{align}
Then $x_0\in S$ is a minimizer of $f$ if and only if it is an efficient element to $\psi$, i.e. $\of{\cb{\psi(x_0)}+(-C)}\cap \bigcup\limits_{x\in S}\psi(x)\subseteq \cb{\psi(x_0)}+C$. A set $M\subseteq X$ is an infimizer if and only if the following domination property is satisfied
\[
\bigcup\limits_{x\in X}f(x)\subseteq\cl\co \bigcup\limits_{m\in M}f(m).
\]

The next result provides some characterizations of minimizers via scalarizations.

\begin{proposition}\label{prop:minimal_iff}
Let $f:X\to\G^\triup$ and $x_0\in\dom f$. Then the following are equivalent
\begin{enumerate}[(a)]
\item 
$
f(x_0)\in\Min f\sqb{X};
$
\item 
$
f(x)\neq f(x_0)\;\Rightarrow\; \exists z^*\in\C:\,  \vp_{f,z^*}(x_0)< \vp_{f,z^*}(x); 
$
\item 
$
f(x)\neq f(x_0)\;\Rightarrow\; \exists z^*\in\C:\,  \vp_{f,z^*}(x)\neq-\infty\,\wedge\,\vp_{f,z^*}(x_0)\idif \vp_{f,z^*}(x)<0; 
$
\item 
$
f(x)\neq f(x_0)\;\Rightarrow\; \exists z^*\in\C:\,  0<\vp_{f,z^*}(x)\idif \vp_{f,z^*}(x_0); 
$
\item
$
f(x)\neq f(x_0)\;\Rightarrow\; 0\notin f(x)\idif f(x_0).
$
\end{enumerate}
\end{proposition}
\proof
Equivalences from (a) through (d) are immediate and by Lemma \ref{lem:scal_of_difference}, (d) and (e) are equivalent.
\pend

\begin{definition}[Solution]\label{def:Infimizer_Minimizer_Solution}\cite{HeydeLoehne11}
Let $f:X\to \G^\triup$. An infimizer of $f$ consisting of only minimizers is called a {\em solution} of the optimization problem \eqref{eq:optim_problem}.
\end{definition}

\begin{example}
Let $f:\R\to\G(\R^2,\R^2_+)$ be given as
\[
f(x)=\cb{(-x,-x)}\oplus R^2_+.
\]
Then $\N\subseteq \R$  as well as any interval $\of{x,+\infty}\subseteq\R$ are infimizers of $f$. However, $\Min f\sqb{\R}=\emptyset$. Hence no solution of $f$ exists.
\end{example}

In \cite{HamelSchrage13PJO} the concept of $z^*$--minimizers was introduced, defining $x_0\in X$ as a $z^*$--minimizer of $f:X\to\G^\triup$ if and only if $x_0$ is a minimizer of $\vp_{f,z^*}:X\to\OLR$. In fact, this concept is independent from the one we are investigating. The following Example  \ref{ex:minimizer_z*Minimizer}(a) due to F. Heyde proves that a solution in the sense of Definition \ref{def:Infimizer_Minimizer_Solution} does not need to be a $z^*$--solution, while Example \ref{ex:minimizer_z*Minimizer}(b) provides a counterexample to the reverse implication.

\begin{example}\label{ex:minimizer_z*Minimizer}
\begin{enumerate}[(a)]
\item
Let $X=Z=\R^2$ and $C=\R^2_+$. The (closed and convex) function $f:X\to\G^\triup$ is defined as follows
\[
f(x)=\begin{cases}
\cb{z\in -x_1+x_2\leq z_1,\, -x_1-x_2\leq z_2,\, x_1\leq z_1+z_2} ,\text{ if } 0\leq x_1;\\
\emptyset,\text{ else.}
\end{cases}
\]
Then each $x_0\in \dom f$ is minimal and $M=\cb{x\in X\st 0<x_1,x_2}$ is a solution of \eqref{eq:optim_problem}, while no $x\in M$ is a $z^*$--solution for any $z^*\in\C$.
\item
Let $X=\R$, $Z=\R^2$ and $C=\R^2_+$. The (closed and convex) function $f:X\to\G^\triup$ is defined as follows
\begin{align*}
f(x)=\begin{cases}
\cb{z\in Z\st \frac{1}{z_1}\leq z_2},\text{ if } 0=x;\\
\cb{z\in Z\st 0\leq z_1,z_2},\text{ if } 1=x;\\
x f(1)\oplus (1-x)f(0),\text{ if } 0\leq x\leq 1;\\
\emptyset,\text{ else.}
\end{cases}
\end{align*}
Then each $x_0\in \dom f$ is  $z^*$--minimal with respect to $z^*\in\cb{(0,-1)^T,(-1,0)^T}$, but the only minimizer of $f$ is $x=1$ and $M=\cb{1}$ is the only solution of \eqref{eq:optim_problem}.
\end{enumerate}
\end{example}


\subsection{Directional derivatives}

The notions of variational inequalities related to an optimization problem involves the concept of directional derivatives.

We apply the following definition to  convex functions $f:X\to\G^\triup$ which  extends the concept of (lower) Dini derivatives to functions mapping to any $\inf$--residuated image space.

We stress that this approach allows to extend the classical Dini derivative for scalar-valued functions to extended real-valued functions (see e.g. \cite{HamelSchrage13PJO, Diss}), as discussed in Example \ref{ex:scalar_dirder} below.

\begin{definition}\label{def:Dirder}
Let $f:X\to\G^\triup$ be convex, $x, u\in X$, then the directional derivative of $f$ at $x$ along direction $u$ is defined as
\begin{align*}
f'(x,u)=\liminf\limits_{t\downarrow 0}\frac{1}{t}\of{f(x+tu)\idif f(x)}
=\bigcap\limits_{0<t_0}\cl\co\bigcup\limits_{t\in\of{0,t_0}}\frac{1}{t}\of{f(x+tu)\idif f(x)}.
\end{align*}
\end{definition}

For convex (set-valued) functions, the differential quotient is monotone.

\begin{proposition}\label{prop:diff_quot_mon}
Let $f:X\to\G^\triup$ be convex, $x_0\in X$ and $g:\of{0,+\infty}\to\G^\triup$ be given  by 
$g(t)=\frac{1}{t}\of{f(x+tu)\idif f(x)}$.
Then for all $0<s\leq t$ it holds $g(s)\lel g(t)$.
\end{proposition}
\proof
Let $z_t\in g(t)$ and $0<s< t$ be satisfied, then there exists an $r\in\of{0,1}$ such that $s=rt$ and 
$f(x+su)\lel (1-r)f(x)\oplus r f(x+tu)$. Thus,
\begin{align*}
f(x+su)\idif f(x)\lel r (f(x+tu)\idif f(x)),
\end{align*}
which in turn implies that
\begin{align*}
\frac{1}{s}\of{f(x+su)\idif f(x)}\lel \frac{r}{rt} (f(x+tu)\idif f(x)),
\end{align*}
as desired.
\pend

The following result extends a well known property of Dini derivatives for convex single-valued functions.

\begin{proposition}\label{prop:f'_pos_hom}\label{prop:f'_conv}
Let $f:X\to\G^\triup$ be convex, $x\in\dom f$ and $u\in X$. Then
\[
f'(x,u)=\inf\limits_{0<t}\frac{1}{t}\of{f(x+tu)\idif f(x)},
\]
$f'(x,0)=0^+f(x)$ and the function $u\mapsto f'(x,u)$ is sublinear as a function from $X$ to $\G(Z,0^+f(x))$.
\end{proposition}
\proof
The first statement comes directly from Proposition \ref{prop:diff_quot_mon}.

For all $x\in X$, $f'(x,0)=\inf\frac{1}{t}\of{f(x)\idif f(x)}$ and thus
\[
f'(x,0)=\begin{cases}
0^+f(x)&\text{, if }x\in \dom f;\\
Z&\text{, elsewhere.}
\end{cases}
\]

By definition, for all $0<s$, $u\in X$ it holds
\[
f'(x,su)=s\cdot\inf\limits_{0<t}\frac{1}{st}\of{f(x+tsu)\idif f(x)}=s f'(x,u).
\]

Let $x,u_1,u_2\in X$ and $s\in\of{0,1}$ be assumed. By Proposition \ref{prop:diff_quot_mon} the differential quotient is decreasing, so for all $0<t_0$ it holds
\begin{align*}
f'(x,su_1+(1-s)u_2)
&=\inf\limits_{0<t\leq t_0}\frac{1}{t}\of{f(s(x+tu_1)+(1-s)(x+tu_2))\idif f(x)}.
\end{align*}
Convexity and Lemma \ref{lem:calc_of_diff} $(b)$  imply
\begin{align*}
f'(x,su_1+(1-s)u_2)
&\lel\inf\limits_{0<t\leq t_0}\frac{1}{t}\of{s\of{f(x+tu_1)\idif f(x)}\oplus (1-s)\of{f(x+tu_2)\idif f(x)}}.
\end{align*}
Since $\G^\triup$ is $\inf$--residuated and by Proposition \ref{prop:diff_quot_mon},
\begin{align*}
f'(x,su_1+(1-s)u_2)
&\lel \frac{1}{t_0}\of{s\of{f(x+t_0u_1)\idif f(x)}}\oplus (1-s)\inf\limits_{0<t\leq t_0}\frac{1}{t}\of{f(x+tu_2)\idif f(x)}\\
&= s\frac{1}{t_0}\of{\of{f(x+t_0u_1)\idif f(x)}}\oplus (1-s)f'(x,u_2).
\end{align*}
But, as this is true for all $0<t_0$ and $\G^\triup$ is $\inf$--residuated, 
\begin{align*}
f'(x,su_1+(1-s)u_2)
&\lel sf'(x,u_1)\oplus (1-s)f'(x,u_2)
\end{align*}
is satisfied.
\pend

\begin{remark}
Since the differential quotients $\frac{1}{t}\of{f(x+tu)\idif f(x)}$ of a convex function $f:X\to\G^\triup$ form a decreasing net of convex sets, their union is convex. Therefore in this case the following equation holds true.
\[
f'(x,u)=\cl\co\bigcup\limits_{t>0}\frac{1}{t}\of{f(x+tu)\idif f(x)}=\cl\bigcup\limits_{t>0}\frac{1}{t}\of{f(x+tu)\idif f(x)}
\]
\end{remark} 

\begin{remark}
Let $f:X\to\G^\triup$ be convex, $x_0\in\dom f$ and $x\in X$. 

If $f'(x_0,x-x_0)\neq\emptyset$, then $\sqb{0,t_0}\subseteq \dom f_{x_0,x}$ is true for some $t_0\in\of{0,1}$ and for all $t\in\of{0,t_0}$ it holds 
\begin{equation*}
0^+f'(x_0,x-x_0)\lel 0^+f(x_t)\lel 0^+f(x_0).
\end{equation*}
Indeed, as $f$ is convex, $0^+f(x_t)$ is constant on the set $\of{0,t_0}$ and $0^+f(x_t)\lel 0^+f(x_0)$.
Also,
\[
f'(x_0,x-x_0)\lel \frac{1}{t}\of{f(x_t)\idif f(x_0)}
\]
and both sets are convex, hence $0^+f'(x_0,x-x_0)\lel 0^+f(x_t)$ by Lemma \ref{lem:Rec_A2}.
\end{remark}

\begin{example}\label{ex:scalar_dirder}
Let $\vp:X\to\OLR$ be convex, $f:X\to\G(\R,\R_+)$ its epigraphical extension as defined in \eqref{eq:sv_extension}. 
If $\vp:X\to\OLR$ is proper, $x\in\dom \vp$, then $f'(x,u)$ coincides with the upper Dedekind cut of the classic directional derivative of $\vp$, while in general,
\begin{align*}
f'(x,u)=\of{\inf\limits_{0<t}\frac{1}{t}\of{\vp(x+tu)\idif \vp(x)}}+\R_+.
\end{align*} 
Especially, if $\vp(x)=+\infty$, then $f'(x,u)=\R$ for all $u\in X$, while if $x\in\dom \vp$ and $\vp(x)=-\infty$, then a careful case study provides
\begin{align*}
f'(x,u)=\begin{cases}
		\R, &\text{ if } u\in\cone\of{\dom \vp-\cb{x}}\\
		\emptyset, &\text{ else.}
\end{cases}
\end{align*}
Therefore
\begin{align*}
\vp'(x,u)=\inf\limits_{0<t}\frac{1}{t}\of{\vp(x+tu)\idif \vp(x)}
\end{align*}
for all $x,u\in X$ provides an extension of Dini derivatives to the case where $\vp$ is improper or $x\notin\dom \vp$.
\end{example}

\begin{remark}\label{rem:dom_of_f'}
Let $f:X\to\G^\triup$ be convex.
It is easy to see that if $x\notin\dom f$, then $f'(x,u)=Z$ and $\vp'_{f,z^*}(x,u)=-\infty$ are satisfied for all $u\in X$ and all $z^*\in\C$.

On the other hand, if $x\in\dom f$, then 
$\dom \vp'_{f,z^*}(x,\cdot)=\cone\cb{\dom f+\cb{-x}}\cup\cb{0}$ is true for all $z^*\in\C$ and the derivative is sublinear. Hence, $\vp'_{f,z^*}(x,u)=-\infty$ implies either $\vp_{f,z^*}(x)=-\infty$, or $\vp'_{f,z^*}(x,-u)=+\infty$.

Especially, $\dom f'(x,\cdot)\subseteq \dom \vp'_{f,z^*}(x,\cdot)$ is always satisfied. Hence if $\vp_{f,z^*}(x)\in\R$, then either $x-tu\notin\dom f$ for all $0<t$, or $-\infty<\vp'_{f,z^*}(x,u)\leq \vp_{f'(x,\cdot),z^*}(u)$.

\end{remark}

If for some $z^*\in \C$ it holds $f(x)=f_{z^*}(x)$ for all $x\in X$ and $f$ is convex, then the scalarization of the derivative  is equal to the  derivative of the scalarization,
$\vp_{f'_{z^*}(x,\cdot),z^*}(u)=\vp'_{f,z^*}(x,u)$ for all $x,u\in X$. However, in general only the following inequality can be proven
\[
\forall z^*\in\C,\;\forall x,u\in X:\quad \vp'_{f,z^*}(x,u)\leq \vp_{f'(x,\cdot),z^*}(u).
\]  

\begin{example}
Let $f:\R\to\G(\R,\cb{0})$ be defined as $f(x)=\sqb{-\sqrt{1-x^2}, \sqrt{1-x^2}}$, whenever $x\in\sqb{-1,1}$ and $f(x)=\emptyset$, else.
Then $f(0)+\cb{z}\nsubseteq f(t)$ for any $t\neq 0$, so $f'(0,u)=\emptyset$.
On the other hand, $\vp_{f,s}(x)=-|s|\cdot \sqrt{1-x^2}$ for all $s\neq 0$ and thus $\vp'_{f,s}(x,u)=-|s|\cdot \frac{x}{\sqrt{1-x^2}}\cdot u$ for all $x\in\of{-1,1}$, especially $\vp'_{f,s}(0,u)=0$ for all $s\neq 0$.
Hence, 
\begin{align*}
\emptyset=f'(0,u)\subsetneq \bigcap\limits_{z^*\in(\cb{0})^-\setminus\cb{0}}f'_{z^*}(0,u)=\cb{0}
\end{align*}
\end{example}

\begin{proposition}\label{prop:vp'fLeqf'}
Let $f:X\to\G^\triup$ be convex and $x,u\in X$. Then 
\begin{align*}
&\bigcap\limits_{z^*\in\C}f'_{z^*}(x,u)\lel f'(x,u);\\
&\forall z^*\in \C:\quad \vp'_{f,z^*}(x,u)\leq \vp_{f'(x,\cdot),z^*}(u) .
\end{align*}
\end{proposition}
\proof
By definition and Lemmas \ref{lem:scal_of_difference} and \ref{prop:scal_of_infimum},
\begin{align*}
f'(x,u)
	&=\cl\co\bigcup\limits_{0<t}\bigcap\limits_{z^*\in\C}
	\cb{z\in Z\st \frac{1}{t}\of{\vp_{f,z^*}(x+tu)\idif \vp_{f,z^*}(x)}\leq -z^*(z)}\\
	&\subseteq\bigcap\limits_{z^*\in\C}\cl\co\bigcup\limits_{0<t}
	\cb{z\in Z\st \frac{1}{t}\of{\vp_{f,z^*}(x+tu)\idif \vp_{f,z^*}(x)}\leq -z^*(z)}\\
	&=\bigcap\limits_{z^*\in\C}
	\cb{z\in Z\st \inf\limits_{0<t}\frac{1}{t}\of{\vp_{f,z^*}(x+tu)\idif \vp_{f,z^*}(x)}\leq -z^*(z)}
\end{align*}
hence the inclusion is proven, implying the inequality as well.
\pend

In the sequel, some results require equality in at least one of the inequalities in Proposition \ref{prop:vp'fLeqf'}.
By {\em strong regularity}, we refer to condition
\begin{equation}\label{eq:reg_sharp}\tag{SR}
\forall z^*\in \C:\quad \vp_{f'(x,\cdot),z^*}(u)=\vp'_{f,z^*}(x,u)
\end{equation}
and by {\em weak regularity} to the following condition.
\begin{equation}\label{eq:reg_weak}\tag{WR}
f'(x,u)= \bigcap\limits_{z^*\in\C}f'_{z^*}(x,u)
\end{equation}
Clearly, \eqref{eq:reg_sharp} implies \eqref{eq:reg_weak}.


\section{Main Results}\label{sec:Main_Res}

As our solution concept involves both attainment of the infimum in  a set and minimality of each element in this set, we need suitable inequalities for each of these properties.
Beginning with the infimizer's part, we need to consider that the solution of a variational inequality is usually a singleton, while the infimizer of \eqref{eq:optim_problem} is a set.
However, Proposition \ref{prop:infimizer} allows to characterize an infimizer $M$ by proving $\hat f(0;M)=\inf f\sqb{X}$, or in other words $\cb{0}$ is a single-valued infimizer of the optimization problem
\begin{align*}
\mbox{minimize} \quad \hat f(x;M) \quad \mbox{subject to} \quad  x \in X.
\end{align*}

Given a single-valued convex function $\vp:X\to\OLR$, a solution to a variational inequality of Stampacchia type is a point $x_0\in X$ such that $0\leq \vp'(x_0,x-x_0)$ for all $x\in X$.
According to our setting, a natural extension of this property is given in the following definition.

\begin{definition}\label{def:str_set_Stamp}
Let $f:X\to\G^\triup$ be convex and $x_0\in\dom f$. Then $x_0$ solves the strict set-valued Stampacchia inequality if and only if
\begin{equation}\label{eq:str_set_Stamp}
\tag{$SVI_I$}
\forall x\in X:\quad  0^+f(x_0)\lel f'(x_0,x-x_0).
\end{equation}
\end{definition}

However, it turns out that, in the set--valued case, infimizers (and minimizers) are often characterized more adequately if a scalar type of variational inequalities is considered.

\begin{definition}\label{def:str_scalar_Stamp}
Let $f:X\to\G^\triup$ be convex, $x_0\in\dom f$. Then $x_0$ solves the strict scalarized Stampacchia inequality if and only if
\begin{equation}\label{eq:str_scalar_Stamp}
\tag{$svi_I$}
\forall x\in X,\,\forall z^*\in \C:\quad 
\vp_{f,z^*}(x_0)=-\infty\,\vee\,0\leq \vp'_{f,z^*}(x_0,x-x_0).
\end{equation}
\end{definition}

Scalarized and set-valued variational inequalities are not equivalent without further assumptions. 

\begin{proposition}\label{prop:Stamp_strict}
Let $f:X\to\G^\triup$ be convex, $x_0\in\dom f$. If $x_0$ solves \eqref{eq:str_scalar_Stamp}, then it also solves \eqref{eq:str_set_Stamp}. If additionally the strong regularity condition \eqref{eq:reg_sharp} is satisfied, then the reverse implication is true as well.
\end{proposition}
\proof
By  Proposition \ref{prop:vp'fLeqf'}, \eqref{eq:str_scalar_Stamp} implies
\[
\bigcap\limits_{\substack{z^*\in\C\\ \vp_{(f,z^*)}(x_0)\neq-\infty}}\cb{z\in Z\st 0\leq -z^*(z)} \lel  \bigcap\limits_{z^*\in\C}f'_{z^*}(x_0,x-x_0)\lel f'(x_0,x-x_0).
\]
By \eqref{eq:RecCone_finiteScal} this implies \eqref{eq:str_set_Stamp} as $\dom f=\dom \vp_{(f,z^*)}$ is true for all $z^*\in\C$.
On the other hand, by Proposition \ref{prop:Rec_A},  \eqref{eq:reg_sharp} combined with \eqref{eq:str_set_Stamp} implies  \eqref{eq:str_scalar_Stamp}.
\pend

\begin{theorem}\label{thm:Stamp_strict}
Let $f:X\to\G^\triup$ be convex, $x_0\in\dom f$. Then $x_0$ solves \eqref{eq:str_scalar_Stamp} if and only if $f(x_0)=\inf f\sqb{X}$.
\end{theorem}
\proof
By Proposition \ref{prop:infimizer_iff} $(e)$, $f(x_0)=\inf f\sqb{X}$ is true if and only if 
\begin{equation*}
\forall x\in X,\,\forall z^*\in \C:\quad 
\vp_{f,z^*}(x_0)=-\infty\,\vee\,0\leq \vp_{f,z^*}(x)\idif \vp_{f,z^*}(x_0),
\end{equation*}
which immediately implies \eqref{eq:str_scalar_Stamp}. The opposite implication is true, as convexity of $f$ implies $\vp_{f,z^*}$ is convex.
\pend

\begin{remark}
According to Proposition \ref{prop:Stamp_strict} and Theorem \ref{thm:Stamp_strict}, the set-valued variational inequality \eqref{eq:str_set_Stamp} is a necessary condition for $\cb{x_0}$ to be an infimizer of $f$. Under the regularity condition   \eqref{eq:reg_sharp} it is also a sufficient condition.
\end{remark}
 
Given a single-valued convex function $\vp:X\to\OLR$, a solution to a variational inequality of Minty type is a point $x_0\in X$ such that $\vp'(x,x_0-x)\leq 0$ for all $x\in X$.

\begin{definition}\label{def:str_set_Minty}
Let $f:X\to\G^\triup$ be convex, $x_0\in\dom f$. Then $x_0$ solves the strict set-valued Minty inequality iff
\begin{equation}\label{eq:str_set_Minty}
\tag{$MVI_I$}
\forall x\in X:\quad f'(x,x_0-x)\lel 0^+f(x_0).
\end{equation}
\end{definition}

Equivalently, $x_0$ is a solution to the strict set-valued Minty inequality if and only if
\begin{equation*}
\forall x\in X:\quad 0\in f'(x,x_0-x).
\end{equation*}

The previous definition can be related to the following family of a scalar Minty inequalities.

\begin{definition}\label{def:str_scalar_Minty}
Let $f:X\to\G^\triup$ be convex, $x_0\in\dom f$. Then $x_0$ solves the strict scalarized Minty inequality iff
\begin{equation}\label{eq:str_scalar_Minty}
\tag{$mvi_I$}
\forall x\in X,\,\forall z^*\in \C:\quad \vp'_{f,z^*}(x,x_0-x)\leq 0.
\end{equation}
\end{definition}

\begin{proposition}\label{prop:SV_Sc_Minty}
Let $f:X\to\G^\triup$ be convex, $x_0\in\dom f$. If $x_0$ solves \eqref{eq:str_set_Minty}, then it also solves \eqref{eq:str_scalar_Minty}. If additionally the regularity condition \eqref{eq:reg_weak} is satisfied, the reverse implication holds true.
\end{proposition}
\proof
If $x_0$ solves \eqref{eq:str_set_Minty}, then Proposition \ref{prop:vp'fLeqf'} implies \eqref{eq:str_scalar_Minty}.
On the other hand, assuming \eqref{eq:str_scalar_Minty} and the regularity condition \eqref{eq:reg_weak}, then $0\in f'(x,x_0-x)$ is satisfied for all $x\in X$, in other words \eqref{eq:str_set_Minty}.
\pend

\begin{theorem}
Let $f:X\to\G^\triup$ be convex, $x_0\in\dom f$. Then $f(x_0)=\inf f\sqb{X}$ if and only if $x_0$ solves \eqref{eq:str_set_Minty} and for all $x\in X$ the function $f_{x_0,x}:\sqb{0,1}\to\G^\triup$ is lattice l.s.c. at $0$.

If $x_0$ solves \eqref{eq:str_scalar_Minty} and for all $x\in X$ 
the function $f_{x_0,x}$ is $\of{\C}$--l.s.c. at $0$, then $f(x_0)=\inf f\sqb{X}$.
\end{theorem}
\proof
By Proposition \ref{prop:infimizer_iff} $(d)$, $f(x_0)=\inf f\sqb{X}$ if and only if $0\in f(x_0)\idif f(x)$ for all $x\in X$, hence by the monotonicity of the differential quotient (see Proposition \ref{prop:diff_quot_mon}) 
\[
f'(x,x_0-x)\lel  f(x_0)\idif f(x)\lel 0^+f(x_0)
\]
is satisfied, proving \eqref{eq:str_set_Minty}. When $f(x_0)\lel f(x)$ is assumed,
\[
f(x_0)\lel \bigcap\limits_{t_0\in\of{0,1}}\cl\co\bigcup\limits_{t\in\of{0,t_0}} f_{x_0,x}(t)
\]
is satisfied and hence $f_{x_0,x}$ is lattice l.s.c. at $0$ for all $x\in X$.

On the other hand, \eqref{eq:str_set_Minty} combined with convexity of $f$ implies
\[
\forall x\in X,\forall s,t\in (0,1]:\quad s<t\;\Rightarrow\; f(x_s)\lel f(x_t).
\] 
Hence if $f_{x_0,x}$ is lattice l.s.c. at $0$, then we obtain
\[
\forall x\in X:\quad f(x_0)=\inf f_{x_0,x}\sqb{0,1}\lel f(x)
\]
and $f(x_0)=\inf f\sqb{X}$ is proven.

The proof of the last implication goes along the same lines.
\pend

Recall that if $f_{x_0,x}:\sqb{0,1}\to\G^\triup$ is $\of{\C}$--l.s.c. at $0$ for all $x\in X$, then each such function is also lattice l.s.c. at $0$. In this case, \eqref{eq:str_set_Minty} and \eqref{eq:str_scalar_Minty} are equivalent.

\begin{remark}
The previous results are summarized in the following scheme of relations.

 \begin{center}
      \includegraphics[width=10 cm]{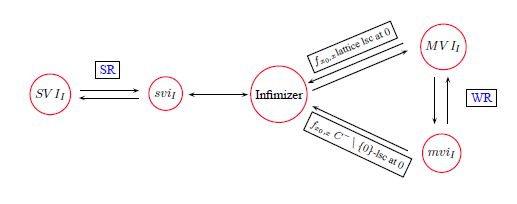}
    \end{center}

\end{remark}

Applying the previous relations and the $\inf$--translation we get a variational characterization of a set $M$ to be an infimizer of $f$.

\begin{corollary}
Let $f:X\to\G^\triup$ be convex, $M\subseteq X$ a set with $M\cap \dom f\neq \emptyset$ and
\[
\hat f(0; M)=\hat f (0;\co M).
\]
Then, $M$ is an infimizer of $f$ if and only if 
 \eqref{eq:str_scalar_Stamp} is satisfied at $0$ for $\hat f(\cdot;\co M)$.
In this case, $\hat f(\cdot;\co M)$ is $\of{\C}$--l.s.c. at $0$ and  \eqref{eq:str_set_Minty}  (and \eqref{eq:str_scalar_Minty})  is satisfied at $0$ for $\hat f(\cdot;\co M)$.

On the other hand, if  \eqref{eq:str_set_Minty} (or \eqref{eq:str_scalar_Minty}) is satisfied at $0$ for $\hat f(\cdot;\co M)$ and one of the conditions $(b)$ or $(c)$ of Lemma \ref{lem:f_M_is convex,lsc} is satisfied, then
$\hat f(\cdot;\co M)$ is $\of{\C}$--l.s.c. at $0$ and $M$ is an infimizer of $f$.
\end{corollary}

In the remainder of this section, we deal with the relation between solutions of variational inequalities and minimizers. 
The variational inequalities of Stampacchia, as well as Minty type are presented both in a set-valued and a scalar(ized) form.

\begin{definition}\label{def:set_Stamp}
Let $f:X\to\G^\triup$ be convex, $x_0\in\dom f$. Then $x_0$ solves the set-valued Stampacchia inequality iff
\begin{equation}\label{eq:set_Stamp2}
\tag{$SVI_M$}
f(x_0)=Z\;\vee\;\forall x\in \dom f:\, f(x)\neq f(x_0)\,\Rightarrow\, 0\notin f'(x_0,x-x_0).
\end{equation}
\end{definition}

\begin{remark}
In \eqref{eq:set_Stamp2}, the condition '$0\notin f'(x_0,x-x_0)$' provides a set-valued version of the property '$\vp'(x_0,x-x_0)\nleq 0$' for scalar convex functions. 
The same inequality could be expressed also by the condition
\begin{equation}\label{eq:set_Stamp}
f(x_0)=Z\;\vee\;\forall x\in \dom f:\, f(x)\neq f(x_0)\,\Rightarrow\,
 f'(x_0,x-x_0)\cap -0^+f(x_0)=\emptyset.
\end{equation}

However, since $\G^\triup$ is not totally ordered, there is a notable difference between these and the condition $f'(x_0,x-x_0)\subset 0^+f(x_0)$.
\end{remark}

\begin{definition}\label{def:scalar_Stamp}
Let $f:X\to\G^\triup$ be convex, $x_0\in\dom f$. Then $x_0$ solves the scalarized Stampacchia inequality, iff
\begin{equation}\label{eq:scalar_Stamp}
\tag{$svi_M$}
f(x_0)=Z\;\vee\;\forall x\in \dom f:\, f(x)\neq f(x_0)\,\Rightarrow\,\exists z^*\in \C:\, 0< \vp'_{f,z^*}(x_0,x-x_0).
\end{equation}
\end{definition}

Property \eqref{eq:scalar_Stamp} also implies
\begin{equation}\label{eq:scalar_Stamp2}
\begin{cases}
\forall x\in\dom f:\quad f(x_0)\neq f(x)\quad \Rightarrow\\ \exists z^*\in \C:\quad
-\infty=\vp_{f,z^*}(x_0)<\vp_{f,z^*}(x)\, \vee\, 0< \vp'_{f,z^*}(x_0,x-x_0).
\end{cases}
\end{equation}
If additionally $f_{x_0,x}:\R\to\OLR$ is $\of{\C}$--l.s.c. at $1$ for all $x\in X$, then \eqref{eq:scalar_Stamp} and \eqref{eq:scalar_Stamp2} are equivalent.

\begin{proposition}\label{prop:SV_andSc_Stamp}
Let $f:X\to\G^\triup$ be convex, $x_0\in\dom f$. If $x_0$ solves \eqref{eq:scalar_Stamp}, then it also solves \eqref{eq:set_Stamp2}. If additionally the regularity condition \eqref{eq:reg_weak} is satisfied, then $x_0$ solves \eqref{eq:set_Stamp2} if and only if it solves \eqref{eq:scalar_Stamp}.
\end{proposition}
\proof
By Proposition \ref{prop:vp'fLeqf'},  \eqref{eq:scalar_Stamp} implies  \eqref{eq:set_Stamp2}.

On the other hand,  \eqref{eq:set_Stamp2} combined with the regularity condition \eqref{eq:reg_weak} implies  \eqref{eq:scalar_Stamp}.

\pend

For the sake of completeness, we quote \cite[Proposition 5.5]{HamelSchrage13PJO}, where it is proven that, if $\dom f\neq \emptyset$, then
\[
f_{z^*}(x)=Z\quad\vee\quad\forall x\in X:\quad 0\leq (\vp_{f, z^*})'(x_0,x-x_0) 
\] 
is equivalent to $f_{z^*}(x_0)=\inf f_{z^*}\sqb{X}$. However, as it has already been shown in Example \ref{ex:minimizer_z*Minimizer}, this concept of optimality is not equivalent to the one investigated in this paper.

\begin{theorem}\label{thm:SVI_Min}
Let $f:X\to\G^\triup$ be convex and $x_0\in \dom f$. If $x_0$ solves \eqref{eq:set_Stamp2} or  \eqref{eq:scalar_Stamp2}, then $f(x_0)\in\Min f\sqb{X}$.
\end{theorem}
\proof
Let $x_0$ be a solution of \eqref{eq:set_Stamp2}, then
\[
f(x)\neq f(x_0)\quad \Rightarrow \quad 0\notin f(x)\idif f(x_0)
\]
is immediate, hence by Proposition \ref{prop:minimal_iff} $(e)$ $x_0$ is a minimizer of $f$.
Assuming \eqref{eq:scalar_Stamp2} is satisfied, then
\[
f(x)\neq f(x_0)\quad \Rightarrow \quad \exists z^*\in \C:\; 0<\vp_{f,z^*}(x)\idif \vp_{f,z^*}(x_0)
\]
is satisfied for all $x\in\dom f$, by Proposition \ref{prop:minimal_iff} $(d)$ implying $f(x_0)\in\Min f\sqb{X}$.
\pend

The reverse implication of Theorem \ref{thm:SVI_Min} is not true, as the following example illustrates.

\begin{example}
Let $\psi:\R\to\R$ be given as $\psi(x)=1$ whenever $-1\leq x\leq 1$ and $\psi(x)=\abs{x}$, elsewhere, $f:X\to\G(\R,\R_+)$ its epigraphical extension. 
The negative dual cone of $C=\R_+$ is the set $\cone\of{\cb{-1}}\cup\cb{0}$ and $\vp_{f,z^*}(x)=-z^*\psi(x)$ for all $z^*\in\C$. Notably, it is sufficient to consider the single scalarization $\vp_{f,-1}:\R\to\OLR$ with $\vp_{f,-1}(x)=\psi(x)$ for all $x\in \R$ and \eqref{eq:reg_sharp} is satisfied.
It holds $f(0)\in\Min f\sqb{X}$, but neither   \eqref{eq:set_Stamp2} nor  \eqref{eq:scalar_Stamp2} are satisfied, as $\psi'(0,-x)=0$ and $f'(0,-x)=\R_+$ holds for all $x\in\R$. 
\end{example}

In a similar way, we approach the Minty type inequalities.

\begin{definition}\label{def:set_Minty}
Let $f:X\to\G^\triup$ be convex, $x_0\in\dom f$. Then $x_0$ solves the set-valued Minty inequality iff
\begin{equation}\label{eq:set_Minty}
\tag{$MVI_M$}
f(x)\neq f(x_0)\quad\Rightarrow\quad  0^+f(x)\not\lel f'(x,x_0-x).
\end{equation}
\end{definition}

Again, \eqref{eq:set_Minty} can be interpreted as the set-valued version of the scalar Minty variational inequality, given by 
\[
\vp(x)\neq \vp(x_0)\quad \Rightarrow \quad 0\nleq\vp'(x,x_0-x),
\]
but it  is significantly different from the condition $0^+f(x)\subset f'(x,x_0-x)$, as $\G^\triup$ is not totally ordered.

\begin{definition}\label{def:scalar_Minty}
Let $f:X\to\G^\triup$ be convex, $x_0\in\dom f$. Then $x_0$ solves the scalarized Minty inequality iff
\begin{equation}\label{eq:scalar_Minty}
\tag{$mvi_M$}
f(x)\neq f(x_0)\quad\Rightarrow\quad\exists z^*\in \C:\quad \vp_{f,z^*}(x)\neq-\infty\,\wedge\, \vp'_{f,z^*}(x,x_0-x)< 0.
\end{equation}
\end{definition}

\begin{proposition}
Let $f:X\to\G^\triup$ be convex, $x_0\in\dom f$. If $x_0$ solves \eqref{eq:set_Minty}, then it also solves \eqref{eq:scalar_Minty}. If additionally the regularity condition \eqref{eq:reg_sharp} is satisfied, then $x_0$ solves \eqref{eq:set_Minty} if and only if it solves \eqref{eq:scalar_Minty}.
\end{proposition}
\proof
If $x_0$ solves \eqref{eq:set_Minty}, then Proposition \ref{prop:vp'fLeqf'} implies \eqref{eq:str_scalar_Minty}.
On the other hand, assuming \eqref{eq:scalar_Minty} and  \eqref{eq:reg_sharp}, then 
for all $x\in X$ with $f(x)\neq f(x_0)$ there exists an element $z\in f'(x,x_0-x)\setminus 0^+f(x)$ (compare Proposition \ref{prop:Rec_A} and Remark \ref{rem:dom_of_f'}), in other words \eqref{eq:set_Minty} is satisfied.
\pend

\begin{proposition}
Let $f:X\to\G^\triup$ be convex and $x_0\in \dom f$.
Then $x_0$ solves \eqref{eq:scalar_Minty} if and only if for all $x\in \dom f$
\begin{equation}\label{eq:monotone_if_Minty}
f(x)\neq f(x_0)\quad \Rightarrow\quad \inf f_{x_0,x}\of{0,1}\lel f(x)\,\wedge\, \inf f_{x_0,x}\of{0,1}\neq f(x) .
\end{equation}
\end{proposition}
\proof
Let $x_0$ be a solution of  \eqref{eq:scalar_Minty}. This is equivalent to state that for each $x\in \dom f$ with $f(x)\neq f(x_0)$ there exists an element $z^*\in\C$ and $t\in\of{0,1}$ such that $\vp_{f,z^*}(x_{t})\idif \vp_{f,z^*}(x)<0$ and $\vp_{f,z^*}(x)\neq-\infty$, or equivalently $\vp_{f,z^*}(x_{t})< \vp_{f,z^*}(x)$.

In this case, \eqref{eq:monotone_if_Minty} is immediate, as 
\[
\inf f_{x_0,x}\of{0,1}\lel\bigcap\limits_{t_0\in\of{0,1}}\cl\bigcup\limits_{t\in\of{t_0,1}} f_{x_0,x}(t)\lel f(x)
\]
by convexity
and $\inf f_{x_0,x}\of{0,1}\lel f(x_t)$, hence strict inclusion is satisfied.

On the other hand,  \eqref{eq:monotone_if_Minty} implies that, if $f(x)\neq f(x_0)$, then there exists $t\in\of{0,1}$ and $z^*\in \C$ such that
$\vp_{f,z^*}(x_{t})< \vp_{f,z^*}(x)$. Hence $\vp_{f,z^*}(x)\neq -\infty$ and
$\vp'_{f,z^*}(x,x_0-x)<0$ are satisfied, as the scalarization $\vp_{f,z^*}:X\to\OLR$ is convex.
\pend

\begin{theorem}\label{thm:finite_Minty}
Let $f:X\to\G^\triup$ be convex and $x_0\in \dom f$.
If $f(x_0)\in\Min f\sqb{X}$, then $x_0$ solves \eqref{eq:scalar_Minty}.
If $x_0$ solves 
\begin{equation}\label{eq:scalar_Minty_finite}
f(x)\neq f(x_0)\quad\Rightarrow\quad\exists z^*\in M^*:\quad \vp_{f,z^*}(x)\neq-\infty\,\wedge\, \vp'_{f,z^*}(x,x_0-x)< 0
\end{equation}
where  $M^*\subseteq \C$ is a finite set.
If additionally $f_{x_0,x}$ is $M^*$--l.s.c. at $0$, then $f(x_0)\in\Min f\sqb{X}$.
\end{theorem}
\proof
Let $f(x_0)\in\Min f\sqb{X}$ be assumed, then by Proposition \ref{prop:minimal_iff} $(c)$
\[
f(x)\neq f(x_0)\quad \Rightarrow\quad \exists z^*\in\C:\; \vp_{f,z^*}(x)\neq-\infty\,\wedge\, \vp_{f,z^*}(x_0)\idif \vp_{f,z^*}(x)<0.
\]
As the differential quotient is decreasing, this implies \eqref{eq:scalar_Minty}.

On the other hand, let \eqref{eq:scalar_Minty_finite} be satisfied and let $\of{\vp_{f,z^*}}_{x,x_0}:\sqb{0,1}\to\OLR$ be l.s.c. at $0$ for all $z^*\in M^*$. Then $f(x)\neq f(x_0)$ and convexity and lower semicontinuity of the scalarizations imply that there exist $z^*\in M^*$ and $t\in[0,1)$ such that 
\[
\inf\of{\vp_{f,z^*}}_{x_0,x}\sqb{0,1}=\vp_{f,z^*}(x_t)<\vp_{f,z^*}(x).
\]
Now either $f(x_t)=f(x_0)$ and $f(x)\not\lel f(x_0)$, or there exist $t_1\in[0,t)$ and $z^*_1\in M^*\setminus\cb{z^*}$ such that
\[
\inf\of{\vp_{f,z^*_1}}_{x_0,x}\sqb{0,1}=\vp_{f,z^*_1}(x_{t_1})<\vp_{f,z^*_1}(x_t)\leq\vp_{f,z^*_1}(x).
\] 
Especially,
\begin{eqnarray*}
&\vp_{f,z^*}(x_t)=-\infty\,\vee\, 0\leq \vp'_{f,z^*}(x_t,x_0-x)\\
&\vp_{f,z^*_1}(x)\neq-\infty\,\wedge\, \vp'_{f,z^*_1}(x,x_0-x)<0
\end{eqnarray*}
are satisfied.
As $M^*$ is finite, there exists  $t_0\in[0,1)$ such that
\begin{eqnarray*}
\exists z^*_0\in M^*:\; \inf\of{\vp_{f,z^*_0}}_{x_0,x}\sqb{0,1}=\vp_{f,z^*_0}(x_{t_0})<\vp_{f,z^*_0}(x);\\
\forall z^*\in M^*:\;0\leq\vp'_{f,z^*}(x_{t_0},x_0-x)\,\vee\,\vp_{f,z^*}(x_{t_0})=-\infty. 
\end{eqnarray*}
Hence especially $f(x_{t_0})=f(x_0)$ and $f(x)\not\lel f(x_0)$.
\pend

Property \eqref{eq:scalar_Minty_finite} implies \eqref{eq:scalar_Minty}, as the relevant set of directions $M^*$ is a subset of $\C$. The reverse implication does not hold and the finiteness assumption in Theorem \ref{thm:finite_Minty} cannot be relaxed, as the following example shows.

\begin{example}
Define $z^*_i=-\frac{1}{i+1}(1,i)^T\in (R^2_+)^-\bs\cb{0}$ for all $i \in \N = \cb{0,1, 2, \ldots}$. Let $f \colon \R \to \G(\R^2,\R^2_+)$ be defined by
\[
\forall x\in \R:\quad f(x)=\bigcap\limits_{i \in \N}\cb{z\in Z\st  -\psi_{z^*_i}(x)\leq -z^*_i(z)}
\]
where
\[
\psi_{z^*_i}(x)=\left\{\begin{array}{lcl}
-(i+1)\min\cb{1-x,i x}&:&\text{ if $x\in\sqb{0,1}$ and $i\in\N$;}\\
+\infty&: &\text{ elsewhere.}
	\end{array}
			\right.
\]
As $\psi_{z^*_i}:\R\to\OLR$ is convex and l.s.c. for all $i\in\N$,  $f$ is $\of{\C}$--l.s.c. and convex, and it is easy to see that $f(0)=f(1)=\R^2_+$. 

Defining $z_i(x)\in\R^2$ by 
\begin{align*}
\forall i \in \N\bs\cb{0} \colon 
\cb{z_i(x)}=\cb{z\in Z\st z^*_{i-1}(z)=\vp_{z^*_{i-1}}(x)}\cap \cb{z\in Z\st z^*_{i}(z)=\vp_{z^*_{i}}(x)}
\end{align*}
then $f(x)=\co\cb{z_i(x)\st i \in \N\bs\cb{0}} + C$ is true for all $x\in\of{0,1}$. This implies that $\vp_{f,z^*_i}(x)=\psi_{z^*_i}(x)$ is true for all $x\in\sqb{0,1}$ and all $i\in \N$ and therefore $f(x)\supsetneq f(0)$ is satisfied for all $x\in\of{0,1}$ so $0$ is no minimizer of $f$.

On the other hand, for any given $x\in\of{0,1}$, it exists an $i \in \N\bs\cb{0}$ such that $x\in\of{\frac{1}{i+1},1}$, hence
$\vp'_{f,z^*_i}(x,0-x)=-(i+1)<0$ and $-i\leq \vp_{f,z^*_i}(x)\neq-\infty$.
Hence the assumptions of Theorem \ref{thm:finite_Minty} are satisfied for $x_0=0$, replacing the finite set $M^*$ by $\C$, although $0$ is no minimizer of $f$.
\end{example}

\begin{remark}
The previous results can be summarized in the following scheme of relations.
   \begin{center}
        \includegraphics[width=10cm]{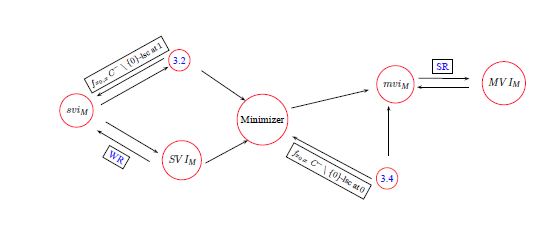}
    \end{center}

\end{remark}


\section{Application to vector optimization}

In this section, we consider a vector-valued function $\psi:S\subseteq X\to Z$ and its epigraphical extension as defined in \eqref{eq:sv_extension}. In the sequel, we refer only to $\dom f=S$, which is the effective domain of $\psi$.

The function $\psi$ is called $C$--convex, when for all $x_1, x_2\in S$ and all $t\in\of{0,1}$ it holds $(1-t)\psi(x_1)+t\psi(x_2)\in \cb{\psi(x_1+t(x_2-x_1))}+C$, or equivalently when 
$\gr f=\epi\psi=\cb{(x,z)\in X\times Z\st z\in \cb{\psi(x)}+C}$
is a convex set, compare \cite[Definition 14.6]{Jahn04}.

\begin{lemma}\label{lem:Deriv_of_Vector}
Let $\psi:S\subseteq X\to Z$ be $C$--convex, $x_0,x\in S$. Then for all $t\in\of{0,1}$ it holds
\[
\frac{1}{t}\of{f(x_0+t(x-x_0))\idif f(x_0)}=\cb{\frac{1}{t}\of{\psi(x_0+t(x-x_0))- \psi(x_0)}}+C.
\]
Moreover $\frac{1}{t}\of{\psi(x_0+t(x-x_0))- \psi(x_0)}$ is decreasing as $t$ converges to $0$
and \eqref{eq:reg_sharp} is satisfied.
\end{lemma}
\proof
By definition, $f(x_t)=\cb{\psi(x_t)}+C$, as  $x_0,x\in S$. Hence 
\[
\forall t\in\of{0,1}:\quad\of{z\in\frac{1}{t}\of{f(x_t)\idif f(x_0)}\quad \Leftrightarrow\quad \psi(x_0)+tz\in \cb{\psi(x_t)}+C},
\]
or equivalently $z\in \cb{\frac{1}{t}\of{\psi(x_0+t(x-x_0))- \psi(x_0)}}+C$.
By Proposition \ref{prop:diff_quot_mon}, the differential quotient is decreasing as $t$ converges to $0$ and by Lemma \ref{prop:scal_of_infimum} 
\[
-\sigma(z^*|f'(x_0,x-x_0))=\inf\cb{-\sigma(z^*| \frac{1}{t}\of{f(x_0+t(x-x_0))\idif f(x_0)})\st 0<t}
\]
for all $z^*\in\C$. But $\vp_{f,z^*}(x)=-z^*\psi(x)$ is satisfied for all $z^*\in\C$ and all $x\in S$, hence
\[
-\sigma(z^*| \frac{1}{t}\of{f(x_0+t(x-x_0))\idif f(x_0)})=-\frac{1}{t}\of{z^*\psi(x_0+t(x-x_0)) - z^*\psi(x_0)},
\]
for all $z^*\in\C$, proving the statement.
\pend

Following the approach in  \cite{CrespiGinRoc2005} we introduce the set of infinite elements $Z_{\infty}=\cb{z_{\infty}\st z\in Z}$.
An element $z_{\infty}$ is the infinite element in direction $z$, in other words
\[
z_{\infty}=\lim\limits_{t\uparrow\infty}tz.
\]
 It holds $z_{\infty}=y_{\infty}$ if and only if $y=\lambda z$ for some $0<\lambda$ and $0_{\infty}=0\in Z$.
For any $z^*\in Z^*$ and $z\in Z$, we define $z^*(z_\infty)=\lim\limits_{t\uparrow +\infty}z^*(tz)$.
Especially,  $z^*(z_{\infty})\in\R$ is satisfied if and only if $z^*(z_{\infty})=z^*(z)=0$.

For a subset $S\subseteq Z$, $S_{\infty}$ denotes the set of all $z_{\infty}\in Z_{\infty}$ with $z\in S\setminus\cb{0}$.
 
The space $\tilde Z=Z\cup Z_\infty$ can be endowed with a topology defined by local bases of neighborhoods as follows.
For any element $z\in Z$, the set $\mathcal U(z)=\mathcal U+\cb{z}$ is a local base of neighborhoods in $\tilde Z$.   
For any element $z\in Z\setminus\cb{0}$, the set
\[
\mathcal U(z_{\infty})=\cb{\of{\cb{tz}+\cone(U+\cb{z})}\cup (U+\cb{z})_{\infty}\st 0<t,\; U\in\mathcal U(z)}
\]
is a local base of neighborhoods of $z_{\infty}$.
Especially, if $K\subseteq Z$ is an open cone with $z\in K$ and $y\in Z$, then $(\cb{y}+K)\cup K_{\infty}$ is a neighborhood of $z_{\infty}$, for details, compare \cite{CrespiGinRoc2005}.

\begin{lemma}
Let $z\in Z$ be given and define $\of{z_{\infty}+C}=\liminf\limits_{t\uparrow\infty}\of{\cb{tz}+C}$,
then $\of{z_{\infty}+C}=\limsup\limits_{t\uparrow\infty}\of{\cb{tz}+C}$ is satisfied.

 If $z\notin -C$, then  $\of{z_{\infty}+C}=\sup\limits_{0<t}\of{\cb{tz}+C}=\emptyset$. Otherwise, 
$\of{z_{\infty}+C}=\inf\limits_{0<t}\of{\cb{tz}+C}$ holds true.

Especially,  $\of{z_{\infty}+C}=C$, if $z\in C\cap -C$ and $\of{z_{\infty}+C}=Z$, if for all $z^*\in \C$ it holds $-z^*(z)<0$. 
\end{lemma}
\proof
By definition,
$\of{z_{\infty}+C}=\bigcap\limits_{0<t_0}\cl\co\bigcup\limits_{t_0\leq t}\of{\cb{tz}+C}.$
Let $z\in -C$, then $\bigcap\limits_{t_0\leq t}\of{\cb{tz}+C}=\cb{t_0z}+C$ and we claim
$\of{\cl\co\bigcup\limits_{0<t_0}\bigcap\limits_{t_0\leq t}\of{\cb{tz}+C}}=\cl\co\bigcup\limits_{0<t_0}\of{\cb{t_0z}+C}$,
or equivalently 
\[
\limsup\limits_{t\uparrow\infty}\of{\cb{tz}+C}=\inf\limits_{t>0}\of{\cb{tz}+C}.
\]
Since
$\inf\limits_{t>0}\of{\cb{tz}+C}\lel \liminf\limits_{t\uparrow\infty}\of{\cb{tz}+C}\lel \limsup\limits_{t\uparrow\infty}\of{\cb{tz}+C}$
always holds true, this implies
\[
\of{z_\infty+C}=\inf\limits_{t>0}\of{\cb{tz}+C}= \limsup\limits_{t\uparrow\infty}\of{\cb{tz}+C}.
\]

On the other hand, let $z\notin -C$ be assumed. Then $0<-z^*(z)$ is satisfied for some $z^*\in \C$. Thus, 
\[
-\sigma(z^*|\cl\co\bigcup\limits_{t_0\leq t}\of{\cb{tz}+C})= -z^*(t_0z)
\]
converges to $+\infty$ as $t_0$ converges to $+\infty$, hence 
$\of{z_{\infty}+C}=\emptyset$.
But, since
\[
\emptyset=\liminf\limits_{t\uparrow\infty}\of{\cb{tz}+C}\lel \limsup\limits_{t\uparrow\infty}\of{\cb{tz}+C}\lel \emptyset
\]
it is proven that 
\[
\of{z_{\infty}+C}=\sup\limits_{t>0}\of{\cb{tz}+C}=\limsup\limits_{t\uparrow\infty}\of{\cb{tz}+C}.
\]
Finally, by Lemma \ref{prop:scal_of_infimum} for $z\in -C$ it holds
\[
\of{z_{\infty}+C}=\bigcap\limits_{z^*\in\C}\cb{y\in Z\st \inf\limits_{0<t}-z^*(tz)\leq -z^*(y)}.
\]
Hence if $z\in C\cap -C$, it is immediate that 
\[
\of{z_{\infty}+C}=\bigcap\limits_{z^*\in\C}\cb{y\in Z\st 0\leq -z^*(y)}=C,
\]
while if for all $z^*\in \C$ it is assumed that $-z^*(z)<0$ holds true, then $\of{z_{\infty}+C}=Z$.

\pend

In \cite{CrespiGinRoc2005} infinite elements play a crucial role to define a Dini directional derivative of $\psi:S\subseteq X\to Z$ at $x_0\in S$ in direction $\of{x-x_0}$ with $x\in S$. The proposed derivative is computed as
\[
\psi'(x_0,x-x_0)=\Limsup\limits_{t\downarrow 0}\cb{\frac{1}{t}\of{\psi(x_0+t(x-x_0))-\psi(x_0)}}\subseteq \tilde Z
\]
where $\Limsup\limits_{t\downarrow 0} \cb{z_{t}}=\cb{\tilde z\in \tilde Z\st \exists \cb{z_{t_i}}_{i\in\N}\subseteq \cb{z_t}_{0<t}, z_{t_n}\to \tilde z}$ 
is the outer Painlev\'e-Kuratowski limit in $\tilde Z$ of a net $\cb{z_{t}}_{t\downarrow 0}\subseteq Z$.

We provide some comparison between the derivative defined in \cite{CrespiGinRoc2005} and our set-valued derivative computed for $\psi^C$.

\begin{lemma}\label{lem:limits_inDeriv_ofVectFct}
Let $\psi:S\subseteq X\to Z$ be $C$--convex,  $f(x)=\psi^C(x)$ for all $x\in X$ and $x_0,x\in S$.

\begin{enumerate}[(a)]
\item If $z\in\psi'(x_0,x-x_0)\cap Z$, then $\cb{z}+C=f'(x_0,x-x_0)$ and for all $z^*\in\C$ it holds $\vp'_{f,z^*}(x_0,x-x_0)=-z^*(z)$;
\item If  $z_{\infty}\in \psi'(x_0,x-x_0)\cap Z_{\infty}$, then $z\in-C$ and $\of{z_{\infty}+C}\subseteq 0^+f'(x_0,x-x_0)$;
\item If $\psi'(x_0,x-x_0)\cap Z\neq\emptyset$ and $z_{\infty}\in\psi'(x_0,x-x_0)\cap Z_{\infty}$, then $z\in C\cap-C$.
\end{enumerate}
\end{lemma}
\proof
\begin{enumerate}[(a)]
\item 
By definition, $z\in\psi'(x_0,x-x_0)\cap Z$ is satisfied if and only if there is a decreasing sequence $\cb{t_i}_{i\in\N}\subseteq \R_+$ such that
$\frac{1}{t_i}\of{\psi(x_0+t_i(x-x_0))-\psi(x_0)}$ converges to $z$ as $i$ converges to $+\infty$. But this implies
\[
\forall z^*\in\C:\quad -z^*(z)\leq \vp'_{f,z^*}(x_0,x-x_0),
\]
hence $\cb{z}+C\supseteq f'(x_0,x-x_0)$. On the other hand, 
\[
z\in\cl\bigcup\limits_{0<t}\of{\cb{\frac{1}{t}\of{\psi(x_0+t(x-x_0))-\psi(x_0)}}+C}=f'(x_0,x-x_0).
\]
\item
Assume to the contrary that $z_\infty\in \psi'(x_0,x-x_0)$ and $z\notin -C$.
Then there exists $U\in\mathcal U$ such that $\cone(U+\cb{z})\cap -C=\emptyset$ and a subsequence 
$z_i=\frac{1}{t_i}\of{\psi(x_0+t_i(x-x_0))-\psi(x_0)}$ with $i\in \N$ 
such that for all $n\in \N$ there exists a $i_0\in \N$ such that for all $i_0\leq i$ it holds $z_i\in \cb{nz}+\cone\of{U+\cb{z}}$, especially $\of{\cb{z_1}+(-C)}\cap \of{\cb{nz}+\cone\of{U+\cb{z}}}\neq\emptyset$ for all $n\in \N$.
However, choosing $n$ sufficiently large, $nz-z_1\in \cone\of{U+\cb{z}}$ is satisfied, implying $\emptyset\neq -C\cap\of{\cb{nz-z_1}+\cone\of{U+\cb{z}}}\subseteq -C\cap\cone\of{U+\cb{z}}=\emptyset$, a contradiction.

\item
Especially by $(a)$, $y\in\psi'(x_0,x-x_0)\cap Z$ is a lower bound of the set 
$$
\cb{\frac{1}{t}\of{\psi(x_0+t(x-x_0))-\psi(x_0)}\st 0<t},
$$
 hence if $z_\infty\in \psi'(x_0,x-x_0)$, then $\forall z^*\in \C:\quad -z^*(y)\leq -z^*(z_{\infty})$, hence by $(b)$ $z\in C\cap-C$. 
\end{enumerate}

\pend

More generally, we remark that taking the limit over a net of singletons and adding the ordering cone does not commute.

\begin{example}
Let $Z=\R^2$, $C=\R^2_+$ be given, $\cb{z_t}_{0<t}\subseteq Z$ a subset of $Z$ with $z_t=(-t,-t^2)$.
Then $\cb{z_t}_{0<t}$ is decreasing as $t$ converges to $+\infty$ and $\Limsup\limits_{t\uparrow+\infty}\cb{z_t}=(0,-1)_{\infty}$.
However, 
\[
\Limsup\limits_{t\uparrow+\infty}\cb{z_t}+C=\cb{z=(z_1,z_2)\in Z\st 0\leq z_1 }\subsetneq \lim\limits_{t\uparrow\infty}\of{\cb{z_t}+C}=Z.
\]
\end{example}

\begin{proposition}\label{prop:tildeZ_compact}\cite{CrespiGinRoc2005}
If $Z$ has finite dimension, then $\tilde Z$ is compact.
\end{proposition}

By Proposition \ref{prop:tildeZ_compact}, if $Z$ has finite dimension, then for a $C$--convex function $\psi:S\subseteq X\to Z$, $x_0,x\in S$ it holds
\[
\emptyset\neq\psi'(x_0,x-x_0)\subseteq Z\cup (-C)_{\infty},
\]
so each element of $\psi'(x_0,x-x_0)$ is either finite (i.e. an element of $Z$), or an element of $(-C)_{\infty}$, (that is an infinite element of $\tilde Z$ which is ''less or equal'' than $0\in Z$ ).

The set of all efficient elements of $\psi\sqb{X}$ is given by
\begin{equation}\label{eq:Eff_Sol}\tag{Eff}
\Eff\psi\sqb{X}=\cb{z\in\psi\sqb{X}\st \forall y\in \psi\sqb{X}:\, z\in \cb{y}+C\,\Rightarrow\, z\in \cb{y}+(-C\cap C)}.
\end{equation}
and $x_0\in\dom f$ is an efficient solution if and only if $\psi(x_0)\in\Eff\psi\sqb{X}$. 
An element $x_0\in\dom f$ is a minimizer of $f$ if and only if it is an efficient solution to $\psi$. Moreover,
\begin{equation}\label{eq:Eff+C=Min}
\bigcup\limits_{f(x)\in\Min f\sqb{X}}f(x)=\Eff\psi\sqb{X}+C
\end{equation}
and a solution to \eqref{eq:optim_problem} exists if and only if
$\cl\co (\Eff\psi\sqb{X}+C)=\cl\co (\psi\sqb{X}+C)$.

In the sequel, we only focus on the characterization of minimizers of $f=\psi^C$ or equivalently efficient solutions of $\psi$.
In this setting, we do not get any new results about infimizer but those already obtained in Section \ref{sec:Main_Res}, as
the $\inf$--translation $\of{\psi^C}\of{\cdot,M}:X\to\G^\triup$ is in general not the epigraphical extension of a vector-valued function.

\begin{corollary}
Let $\psi:S\subseteq X\to Z$ be $C$--convex, $x_0\in S$ and $f(x)=\psi^C(x)$ for all $x\in X$. Then \eqref{eq:set_Stamp2}, \eqref{eq:scalar_Stamp} and \eqref{eq:scalar_Stamp2} are equivalent.  
Especially, if for all $x\in S$ with $\psi(x)\neq \psi(x_0)$ there exists  $z\in Z$ such that $z\in \psi'(x_0,x-x_0)\setminus-C$, then $\psi(x_0)\in\Eff \psi\sqb{X}$.
\end{corollary} 
\proof
The first part of the statement is proven in Proposition \ref{prop:SV_andSc_Stamp}, as by Lemma \ref{lem:Deriv_of_Vector}, \eqref{eq:reg_sharp} and hence especially \eqref{eq:reg_weak} are satisfied. The existence of $z\in Z$ with $z\in \psi'(x_0,x-x_0)\setminus-C$ implies  the existence of a $z^*\in\C$ with $0<\vp'_{f,z^*}(x_0,x-x_0)$, compare Lemma \ref{lem:limits_inDeriv_ofVectFct} $(a)$. Thus  \eqref{eq:scalar_Stamp2} is satisfied, proving the statement.
\pend

\begin{corollary}\label{cor:sv_sc_minty_vv}
Let $\psi:S\subseteq X\to Z$ be $C$--convex, $x_0\in S$ and $f(x)=\psi^C(x)$ for all $x\in X$. Then $x_0$ solves \eqref{eq:set_Minty} if and only if it solves \eqref{eq:scalar_Minty}. Moreover, \eqref{eq:set_Minty} is equivalent to
\begin{equation}\label{eq:sc_Minty_2}
\of{x\in S,\; t\in\of{0,1},\;\psi(x_t)\neq \psi(x_0)}
\;\Rightarrow\; \exists z^*\in\C:\; (-z^*\psi)'(x_t,x_0-x)<0.
\end{equation}
\end{corollary}
\proof
The first part of the statement is true as \eqref{eq:reg_sharp} is guaranteed by Lemma \ref{lem:Deriv_of_Vector} (compare Proposition \ref{prop:SV_Sc_Minty}).
As $f(x)=\psi^C(x)$ for all $x\in X$ is assumed, $\vp_{f,z^*}(x)\neq-\infty$ is always true for all $z^*\in\C$. It is left to prove that \eqref{eq:sc_Minty_2} implies  \eqref{eq:scalar_Minty}.

Let $x\in S$ and $\psi(x_t)\neq \psi(x_0)$ be assumed for some $t\in\of{0,1}$. By convexity of $\vp_{f,z^*}:X\to\OLR$, $(-z^*\psi)'(x_t,x_0-x)<0$ implies $(-z^*\psi)'(x,x_0-x)<0$.
On the other hand, if $\psi(x_t)= \psi(x_0)$ is satisfied for all $t\in\of{0,1}$, then by convexity of the scalarizations \[
-z^*\psi(x_0)=\liminf\limits_{t\downarrow 0}(-z^*\psi(x_t))\leq -z^*\psi(x)
\]
is satisfied for all $z^*\in\C$. Especially, $\psi(x)\neq \psi(x_0)$ implies
\[
\exists z^*\in\C:\quad -z^*(x_0)=-z^*\psi(x_t)<-z^*\psi(x),
\]
hence $\vp'_{f,z^*}(x,x_0-x)=-\infty<0$.
\pend

\begin{remark}
As $(-z^*\psi)'(x,\cdot):X\to\OLR$ is sublinear, if $\psi:S\subseteq X\to Z$ is $C$--convex, $x_0, x\in S$ implies $(-z^*\psi)'(x_t,x_0-x)\in\R$ for all $z^*\in \C$ and all $t\in\of{0,1}$.
In this case, $z_{\infty}\in \psi'(x_t,x_0-x)$ implies $z\in C\cap -C$.

Indeed, under the given assumptions, $-z^*\psi(x_t)\in\R$ is true for all $t\in\of{0,1}$, hence 
\[
0=(-z^*\psi)'(x_t,0)\leq (-z^*\psi)'(x_t,x-x_0)\isum (-z^*\psi)'(x_t,x_0-x)
\]
and $(-z^*\psi)'(x_t,x_0-x)=-\infty$ implies $(-z^*\psi)'(x_t,x-x_0)=+\infty$. But as  $\dom (-z^*\psi)'(x_t,\cdot)=\cone\of{S+\cb{-x_t}}$, this is a contradiction.
By Lemma \ref{lem:limits_inDeriv_ofVectFct} $(b)$, $z_{\infty}\in \psi'(x_t,x_0-x)$ implies $z\in -C$. Assuming $z\notin C$ would imply the existence of a $z^*\in\C$ such that $\psi'(x_t,x_0-x)=-\infty$, a contradiction.
\end{remark}

\begin{proposition}\label{prop:inner_Minty_vv}
Let $\psi:S\subseteq X\to Z$ be $C$--convex, $x_0\in S$ and $f(x)=\psi^C(x)$ for all $x\in X$. If  $x\in S$ and $t\in\of{0,1}$ imply
\[
\psi(x_t)\neq\psi(x_0) \quad\Rightarrow\quad \psi'(x_t,x_0-x)\nsubseteq \of{C\cup C_{\infty}},
\]
 then  $x_0$ solves \eqref{eq:set_Minty} and 
\[
\psi(x_t)\neq\psi(x_0) \quad\Rightarrow\quad \psi'(x_t,x_0-x)\subseteq \of{C\cap -C}_{\infty}\cup\of{Z\setminus C}.
\]
\end{proposition}
\proof
Under the given assumptions,  let $\psi(x_t)\neq\psi(x_0)$. Then $\psi'(x_t,x_0-x)\neq\emptyset$ and especially,
\[
\psi'(x_t,x_0-x)\cap \of{ ((-C)_{\infty}\setminus C_\infty) \cup(Z\setminus C) }\neq\emptyset.
\]
Thus if $z\in \psi'(x_t,x_0-x)\cap (Z\setminus C)$, then there exists an element $z^*\in \C$ satisfying $\vp'_{f,z^*}(x_t,x_0-x)<0$. On the other hand, if $z_\infty \in \psi'(x_t,x_0-x) \cap \of{ ((-C)_{\infty}\setminus C_\infty)}$, then $\vp'_{f,z^*}(x_t,x_0-x)=-\infty$ is satisfied for some $z^*\in\C$, a contradiction. Hence
\[
\emptyset\neq \psi'(x_t,x_0-x)\subseteq ((-C)_{\infty}\cap C_\infty) \cup Z
\]
and thus by assumption
\[
\emptyset\neq \psi'(x_t,x_0-x)\cap \of{Z\setminus C}.
\]
But this implies 
\[
\forall z\in \psi'(x_t,x_0-x)\cap \of{Z\setminus C}:\quad
\emptyset\neq \psi'(x_t,x_0-x)\cap Z\subseteq \cb{z}+(C\cap -C)\subseteq Z\setminus C,
\]
implying the existence of a $z^*\in \C$ satisfying $\vp'_{f,z^*}(x_t,x_0-x)<0$, hence \eqref{eq:scalar_Minty} and therefore \eqref{eq:set_Minty} is satisfied.
\pend

We can prove that under certain assumptions the efficient solutions of a vector valued function are identical with the solutions to the set-valued Minty variational inequality of its epigraphical extension. 

\begin{theorem}\label{thm:finite_Minty_vv}
Let $\psi:S\subseteq X\to Z$ be $C$--convex, $x_0\in S$ and $f(x)=\psi^C(x)$ for all $x\in X$. If $f_{x_0,x}$ is $\of{\C}$--l.s.c. at $0$ for all $x\in X$  and $C$ is polyhedral, then  $x_0$ solves \eqref{eq:set_Minty} if and only if $\psi(x_0)\in\Eff\psi\sqb{X}$.
\end{theorem}
\proof
If $C$ is polyhedral, then so is $C^-$, that is there exists a finite set $M^*=\cb{m_1,...,m_n}\in \C$ such that 
\[
C=\bigcap\limits_{i=1}^{n}\cb{z\in Z\st 0\leq -m^*_i(z)}.
\]
Also, for all $z^*\in \C$, $z^*\in \cone\co M^*$ and for all $z\in Z$ and all $z^*\in \C$, if
\[
z^*=\sum_{i=1}^{n}t_im^*_i,\;0\leq  t_1,...,t_n,
\]
then $-z^*(z)=-\sum_{i=1}^{n}t_im^*_i(z)$.
Let $(-z^*\psi)'(x,x_0-x)<0$ be satisfied for some $z^*=\sum_{i=1}^{n}t_im^*_i\in \C$ and $x_0\in S$.
Then there exists $0<\bar s$ such that (for all $s\in\of{0,\bar s}$)
\[
-z^*(\frac{1}{s}\of{\psi(x+s(x_0-x))-\psi(x)})<0,
\]
hence there exists at least one $i\in\cb{1,...,n}$ such that
\[
-m_i^*(\frac{1}{s}\of{\psi(x_t+s(x_0-x))-\psi(x_t)})<0,
\]
implying $(-m_i^*\psi)'(x,x_0-x)<0$.
In this case, \eqref{eq:scalar_Minty} implies \eqref{eq:scalar_Minty_finite}, thus they are equivalent. Moreover, by Corollary \ref{cor:sv_sc_minty_vv}, \eqref{eq:scalar_Minty} and \eqref{eq:set_Minty} are equivalent. 
As $\psi(x_0)\in\Eff\psi\sqb{X}$ is satisfied if and only if $f(x_0)\in\Min f\sqb{X}$, Theorem \ref{thm:finite_Minty} proves the statement.
\pend

Theorem \ref{thm:finite_Minty_vv} provides as special case the following Minty variational principle for vector-valued functions, which can be found in e.g. \cite{CreGinRoc08,YangYang04}.

\begin{corollary}
Let $Z=\R^m$ and $C=\R^m_+$. Let $\psi:S\subseteq X\to Z$ be $C$--convex, $x_0\in S$ and $f(x)=\psi^C(x)$ for all $x\in X$. If $f_{x_0,x}$ is $\of{\C}$--l.s.c. at $0$ for all $x\in X$, then $\psi(x_0)\in\Eff\psi\sqb{X}$ is satisfied if and only if $x\in S$ and $t\in\of{0,1}$ imply
\[
\psi(x_t)\neq\psi(x_0) \quad\Rightarrow\quad \psi'(x_t,x_0-x)\subseteq Z\setminus C.
\]
Especially in this case, $\psi'(x_t,x_0-x)\subseteq Z$ is single-valued.
\end{corollary}
\proof
By Proposition \ref{prop:tildeZ_compact}, $\psi'(x_t,x_0-x)\neq\emptyset$ is satisfied under the given assumptions and $C$ is polyhedral and pointed, i.e. $C\cap-C=\cb{0}$. Thus 
$
\emptyset\neq\psi'(x_t,x_0-x)\subseteq Z
$
holds true for all $x\in S$ and all $t\in\of{0,1}$ and $\psi'(x_t,x_0-x)$ is single-valued.
Hence, $\psi'(x_t,x_0-x)\subseteq Z\setminus C$ is equivalent to $\psi'(x_t,x_0-x)\nsubseteq  \of{C\cup C_{\infty}}$.
Moreover, under the given assumptions \eqref{eq:set_Minty} is satisfied  (compare Proposition \ref{prop:inner_Minty_vv}). 
By  Theorem \ref{thm:finite_Minty_vv}, \eqref{eq:set_Minty} is equivalent to   $\psi(x_0)\in\Eff\psi\sqb{X}$.

On the other hand, by Corollary \ref{cor:sv_sc_minty_vv}, \eqref{eq:set_Minty} is equivalent to \eqref{eq:sc_Minty_2}, implying
\[
t\in\of{0,1},\, \psi(x_t)\neq \psi(x_0)\quad\Rightarrow\quad \psi'(x_t,x_0-x)\setminus C\neq\emptyset,
\] 
which in turn implies 
\[
t\in\of{0,1},\, \psi(x_t)\neq \psi(x_0)\quad\Rightarrow\quad \psi'(x_t,x_0-x)\subseteq Z\setminus C,
\] 
as proposed.
\pend

%
%

\section{Conclusion}

By means of conlinear spaces we developed a variational inequalities scheme to characterize solutions to set optimization problems. The results proved actually allow to recover results previously proved in vector optimization under convexity assumptions.
It is an open question how far the convexity assumption can be relaxed for set-valued problems.\\

The graphics in the paper summarize the implications proved. Counterexamples are given for the equivalences that do not hold for the formulation presented in the paper.

\section{Acknowledgment}

We thank the two anonymous referee for their valuable comments that helped to improve the manuscript.


\addcontentsline{toc}{section}{Bibliography}

\end{document}